\documentclass[11pt,a4paper,english,reqno]{amsart}
\usepackage{amsmath,amssymb,amsfonts,epsfig,mathrsfs,hyperref}
\usepackage[T1]{fontenc}

\usepackage{color}
\usepackage{array}
\usepackage{amsthm}
\usepackage{amstext}
\usepackage{graphicx}
\usepackage{setspace}

\usepackage[margin=2.5cm]{geometry}
\usepackage{bbm}
\usepackage{color}
\usepackage{enumitem}
\allowdisplaybreaks[4]

\usepackage{amscd,psfrag}
\usepackage[mathscr]{eucal}

\usepackage{slashed}

\makeatletter
\pdfpageheight\paperheight
\pdfpagewidth\paperwidth

\usepackage{comment}

\usepackage{epstopdf}
\usepackage{chngcntr}
\counterwithin{figure}{section}
\usepackage{indentfirst}	

\usepackage[normalem]{ulem}
\theoremstyle{plain}

\newtheorem{definition}{Definition}[section]
\newtheorem{theorem}[definition]{Theorem}
\newtheorem*{theorem*}{Theorem}

\newtheorem{remark}[definition]{Remark}

\newtheorem*{remark*}{Remark}
\newtheorem*{sideremark*}{Side Remark}\newtheorem*{mt*}{Main Theorem}

\newtheorem*{claim*}{Claim}
\newtheorem*{q*}{Question}

\newtheorem{corollary}[definition]{Corollary}
\newtheorem*{corollary*}{Corollary}
\newtheorem*{proposition*}{Proposition}

\newtheorem{proposition}[definition]{Proposition}

\newcommand{\R}{\mathbb{R}}

\newcommand{\na}{\nabla}
\newcommand{\dd}{{\rm d}}
\newcommand{\p}{\partial}
\newcommand{\e}{\epsilon}
\newcommand{\emb}{\hookrightarrow}
\newcommand{\weak}{\rightharpoonup}
\newcommand{\map}{\rightarrow}
\newcommand{\G}{\Gamma}
\newcommand{\M}{\mathcal{M}}
\newcommand{\K}{\mathcal{K}}

\newcommand{\two}{{\rm II}}

\newcommand{\loc}{{\rm loc}}
\newcommand{\lie}{{\mathfrak{g}}}
\newcommand{\conn}{{\mathfrak{A}}}
\newcommand{\T}{{\mathcal{T}}}

\newcommand{\curv}{{\mathscr{F}}}
\newcommand{\nnn}{{\mathfrak{n}}}
\newcommand{\ttt}{{\mathfrak{t}}}
\newcommand{\bra}{\left\langle\left\langle}
\newcommand{\ket}{\right\rangle\right\rangle}

\newcommand{\hato}{{\hat{\Omega}}}

\newcommand{\U}{{\mathcal{U}}}
\newcommand{\V}{{\mathcal{V}}}
\newcommand{\tilo}{{\widetilde{\Omega}}}
\newcommand{\so}{{\mathfrak{so}}}

\newcommand{\harm}{{\mathscr{H}}}

\newcommand{\stwo}{{\mathbf{S}^2}}

\allowdisplaybreaks[4]

\def\XXint#1#2#3{{\setbox0=\hbox{$#1{#2#3}{\int}$ }
\vcenter{\hbox{$#2#3$ }}\kern-.6\wd0}}

\numberwithin{equation}{section}
\numberwithin{figure}{section}

\usepackage{chngcntr}
\counterwithout{equation}{section}

\title{Smoothability of $L^p$-connections on bundles and isometric immersions with $W^{2,p}$-regularity}

\author{Siran Li}

\address{Siran Li: School of Mathematical Sciences and CMA-Shanghai, Shanghai Jiao Tong University, No.~6 Science Buildings, 800 Dongchuan Road, Minhang District, Shanghai, China (200240)}

\email{\texttt{siran.li@sjtu.edu.cn}}



\date{\today}

\begin{document}

\maketitle

\begin{center}
\emph{Dedicated to Constantine M. Dafermos on the occasion of his $85^{\text{th}}$ birthday,\\ with gratitude and admiration}
\end{center}

\begin{abstract}
We are concerned with two interrelated problems: smoothability of connection 1-forms with low regularity on  bundles with prescribed smooth curvature 2-forms, and existence of isometric immersions with low regularity. We first 
show that if $\Omega$ is an $L^p$-connection $1$-form on a vector bundle over a closed Riemannian $n$-manifold $\M$ with small $L^p$-norm ($p>n$) and smooth curvature $2$-form $\curv$, then $\Omega$ can be approximated in the $L^p_\loc$-topology by smooth connections of the same curvature (not necessarily gauge equivalent).  Our proof, adapted from S. Mardare's work \cite{m05} on the fundamental theory of surfaces with $L^p$-second fundamental form, is elementary in nature and uses only Hodge decomposition and fixed point theorems. This result is then applied to the study of isometric immersions of  Riemannian manifolds with low regularity. We revisit the proof for the  existence of $W^{2,p}$-isometric immersion $\M^n \emb \R^{n+k}$ with arbitrary $n$ and $k$ given weak solutions to the Gauss--Codazzi--Ricci equations, aiming at elucidating some global \emph{vs.} local issues. We also give a characterisation for Riemannian metrics that admit $W^{2,p}$- but no $C^\infty$-isometric immersions. 
\end{abstract}

\section{Introduction}\label{sec: Intro}

This note is concerned with two classical problems in geometric analysis: characterisation of the moduli space of connection 1-forms with low regularity on vector bundles/principal G-bundles given prescribed smooth curvature 2-forms, and the existence of isometric immersions with low regularity. Our goal is two-fold: the first is to demonstrate that the nonlinear smoothing techniques developed by S. Mardare (\cite{m05, m07}, motivated in turn by the isometric immersion problem) can be adapted to investigate the moduli space of connections on bundles over an $n$-dimensional manifold $\M$ with only $L^p$-regularity for $p>n$. This goes further below the regularity threshold in Uhlenbeck's classical works \cite{u', u}, in which $W^{1,p}$-connections with $p>\frac{n}{2}$ are considered. The second goal is to study the existence of isometric immersions of an $n$-dimensional manifold $\M$ with $L^p$-extrinsic geometry ($p>n$) into Euclidean spaces. In particular, we clarify some existence results obtained recently in \cite{m05, m07, sz, cl}, and also analyse the examples of $W^{2,p}$- but non-$C^\infty$-isometric immersions.

Besides its interest for geometrically oriented readers, this note is also strongly motivated by applied problems in various fields of applied mathematics and physics --- in particular, in nonlinear elasticity and mathematical relativity. In elasticity theory, one  major objective is to determine the deformation undergone by an elastic body in response to external forces and boundary conditions. The elastic body and deformation are modelled respectively by a manifold $\M^n$ and a mapping $\iota: {\M}^n \to \R^{n+k}$. The \emph{intrinsic approach to nonlinear elasticity}, pioneered by Antman \cite{antman} and Ciarlet \cite{cgm}, recasts the problems concerning the deformation $\iota$ into those concerning the right Cauchy--Green tensor $g=\iota^\#\delta$, \emph{i.e.}, the pullback of the Euclidean metric $\delta$ on $\R^{n+k}$ under $\iota$. Thus, the existence of $\iota$, in geometric terms, amount to the existence of an isometric immersion. On the other hand, the issue of smoothability of $L^\infty$-connections on tangent bundles of arbitrary (pseudo-)Riemannian manifolds with $L^\infty$-curvature has been studied systematically by Reintjes--Temple in \cite{rt1, rt2, rt3, rt4}, establishing that the Lorentzian metrics of shock wave solutions of the Einstein-Euler equations are non-singular.

\subsection{Smoothability of $L^p$-connections}
We first prove a result concerning the structure of the space of $L^p$-affine connections on vector bundles over an $n$-dimensional Riemannian manifold $\M$. Assume throughout $p>n$. Our set up follows closely the seminal works \cite{u, u'} by Uhlenbeck (see also the exposition \cite{w} by Wehrheim). However, as aforementioned, in our note less stringent regularity assumptions shall be imposed.

Throughout this note, our use of Sobolev norms and other geometric-analytic notations on Riemannian manifolds and bundles are standard. See, for example, \cite{d, h, sch} for references. We collect some relevant elementary materials on differential geometry in Appendix~\ref{appendix}.


As in \cite{u, u'}, let $\eta$ be a vector bundle with structure group $G$ over a Riemannian $n$-dimensional manifold $\M$. Let $\conn(\M)$ be the affine space of smooth connections on $\M$:
\begin{equation*}
\conn(\M)= \hato + C^\infty\left(\M; {\rm Ad}\,\eta  \otimes \bigwedge^1T^*\M\right),
\end{equation*}
where $\hato$ is a base connection. Locally on a trivialised chart $\U \subset \M$ the adjoint bundle $${\rm Ad}\,\eta\big|_\U \cong \U \times \lie,$$ where $\lie$ is the Lie algebra of $G$. In other words, $\conn(\M)$ is the affine space modelled over $\lie$-valued differential $1$-forms. We also denote by $\conn^{p}_0(\M)$ the space of connection $1$-forms with $L^{p}$-regularity: $\conn^{p}_0(\M) = \hato + L^p\left(\M; {\rm Ad}\,\eta  \otimes \bigwedge^1T^*\M\right)$. More generally, for $p \geq 1$, $k \geq 0$, write
\begin{equation*}
\conn^{p}_k(\M) = \hato + W^{k,p}\left(\M; {\rm Ad}\,\eta  \otimes \bigwedge^1T^*\M\right).
\end{equation*}
This definition is independent of the choice of $\hato$, which shall be fixed from now on.

The curvature $2$-form of a connection $\Omega \in \conn(\M)$ is defined as
\begin{equation}\label{curvature, def}
F_\Omega = d\Omega + [\Omega \wedge \Omega] \in C^\infty\left(\M; {\rm Ad}\,\eta  \otimes \bigwedge^2T^*\M\right).
\end{equation}
The notation $[\Omega \wedge \Omega]$ designates the intertwining of  wedge product in the $T^*\M$ factor and  commutator in the ${\rm Ad}\,\eta $ (or $\lie$) factor. Eq.~\eqref{curvature, def}  extends naturally to $\Omega \in \conn^{p}_0(\M)$ for any $p \geq 2$. In this case
$F_\Omega \in \left[W^{-1,p} + L^{p/2}\right]\left(\M; {\rm Ad}\,\eta  \otimes \bigwedge^2T^*\M\right)$,  well-defined in the sense of distributions. 

Our first main result is that the moduli space of connection 1-forms with uniformly small $L^p$-norm ($p>n$) and prescribed smooth curvature $2$-form is smoothable. That is, each $L^p$-connection $\Omega$ has a family of smooth connection 1-forms $\{\Omega^\e\}$ that converges in the strong $L^p$-topology to $\Omega$, such that all of $\Omega^\e$ have the same curvature. Under the smallness condition $\|\Omega\|_{L^{p}(M)} < \kappa_0$, this result does not require any assumption on the topology of $\M$.

\begin{theorem}\label{thm: main}
Let $\eta$ be a vector bundle with structure group $G$ over an $n$-dimensional closed Riemannian manifold $\M$; $p>n$. Let $\curv \in W^{s,q}\left(\M;{\rm Ad}\,\eta \otimes \bigwedge^2 T^*\M\right)$ be a prescribed $2$-form with $q\geq\frac{p}{2}$ and $s \geq 0$ such that $W^{s,q}\subset L^2$ over $\M$. There exists a constant $\kappa_0>0$ depending only on $p$ and $\M$ such that the following space is smoothable:
\begin{align*}
\mathfrak{A}^p_0\left(\M,\kappa_0,\curv\right) :=\bigg\{ \Omega \in  \mathfrak{A}^p_0(\M):\, \|\Omega\|_{L^p(\M)} < \kappa_0 \text{ and } F_\Omega = \curv \bigg\}.
\end{align*}
By this we mean that there exists a family of connections $\{\Omega^\e\} \subset \conn^{1+s}_q(\M)$ such that $\Omega^\e \to \Omega$ in $L^p\left(\M;{\rm Ad}\,\eta\otimes \bigwedge^1 T^*\M\right)$ as $\e \to 0$ and $F_{\Omega^\e}\equiv\curv$ for each $0<\e\leq 1$.
\end{theorem}

\begin{remark}
In Theorem~\ref{thm: main}, the regularity class $W^{s,q}$  of $\curv$ can be replaced by any function space that continuously embeds into $L^{\max\left\{2,\frac{p}{2}\right\}}$ over $\M$, provided that the Calder\'{o}n--Zygmund-type elliptic estimates holds in the scale of such function spaces; \emph{i.e.}, $\Delta^{-1}$ improves the regularity by order 2. The reason will become transparent along the proof.
\end{remark}

By virtue of the Sobolev embedding $W^{1,\frac{p}{2}}(\M^n)\emb L^p(\M^n)$ for $p>n$, the family $\{\Omega^\e\}$ of smooth approximations in the statement of Theorem~\ref{thm: main} actually lies in $\mathfrak{A}^p_0\left(\M,\kappa_0,\curv\right)$. In particular, if the prescribed curvature 2-form $\curv \in C^\infty$ --- \emph{e.g.}, when $\curv \equiv 0$ (the connections are flat) --- then we can take $\{\Omega^\e\}\subset \conn(\M)$. This immediately leads to the following:
\begin{corollary}[Regularity pending uniqueness]
In the setting of Theorem~\ref{thm: main}, if the solution to the equation $d\Omega + [\Omega \wedge \Omega]=\curv$ is unique in  $\mathfrak{A}^p_0\left(\M,\kappa_0,\curv\right)$, then  $\mathfrak{A}^p_0\left(\M,\kappa_0,\curv\right)= \mathfrak{A}\left(\M,\kappa_0,\curv\right).$\label{cor: regularity}
\end{corollary}

Theorem~\ref{thm: main} complements the classical strong compactness theorem of Uhlenbeck \cite{u'} (see also \cite[Chpt. III]{w} for an alternative proof by Salamon), which ascertains that any Yang--Mills connection with $W^{1,p}$-regularity ($p>\frac{n}{2}$) is gauge-equivalent to a $C^\infty$-connection modulo $W^{2,p}$-gauge transforms. Our result pertains to $L^p$-connections ($p>n$) with a prescribed curvature 2-form $\curv$. Here $\curv$ is assumed to be fixed but not necessarily Yang--Mills; \emph{i.e.}, we do not impose any PDE constraint on $\curv$. Theorem~\ref{thm: main} may be compared with Wehrheim \cite{w-new}, which establishes the compactness of the moduli space of flat $L^p$-connections and applies this to deal with the Lagrangian boundary conditions for anti-self-dual instantons.

The crucial difference between our result and the classical theorems in \cite{u, u', w} is the following: we do \emph{not} assert that the smooth approximations $\Omega^\e$  in Theorem~\ref{thm: main} are gauge equivalent to the $L^p$-connection  $\Omega$. That is, it is \emph{not} required that $\Omega^\e = \left(u^\e\right)^\# \Omega = \left(u^\e\right)^{-1}\hato u + \left(u^\e\right)^{-1}\Omega u^\e$ for some $G$-bundle automorphism $u^\e \in {\rm Aut}(\eta)$. Correspondingly, our proof for Theorem~\ref{thm: main} does not involve constructions of Uhlenbeck--Coulomb gauges as in \cite{u}; instead, we only utilise elementary arguments via fixed point theorems. 
Our proof essentially follow S. Mardare \cite{m05, m07} (see also Ciarlet--Gratie--Mardare \cite{cgm}; Chen--Li \cite{cl}) on the solubility of Pfaff and Poincar\'{e} systems and the extension of the fundamental theorem of surface theory to $W^{1,p}_\loc$-metrics with $p>2$. 

Theorem~\ref{thm: main} can be made global without the $L^p$-smallness assumption on the curvature 1-form ``in the Luzin sense'' as follows:
\begin{proposition}\label{propn: luzin}
Let $\eta$ be a vector bundle with structure group $G$ over an $n$-dimensional closed Riemannian manifold $\M$; $p>n$. Let $\curv \in \G\left({\rm Ad}\,\eta \otimes \bigwedge^2 T^*\M\right)$ be a prescribed $2$-form with any regularity. Suppose for arbitrarily small $v_0>0$ the following holds: 
\begin{quote}
On any local chart $\U \subset \M$ with Riemannian volume less than $v_0$, one can find a local connection $1$-form $\Omega_\U \in \conn^p_0(\U)$ such that $d \Omega_\U + \left[\Omega_\U\wedge \Omega_\U\right]=\curv\big|_\U$.
\end{quote}
Then, for any given $\e>0$, there exists a global connection $1$-form $\Omega \in \conn^p_0(\M)$ such that
\begin{align*}
&\text{$d\Omega + [\Omega \wedge \Omega] = \curv$ in the sense of distributions}\\
&\text{on a set $\M_\circ \subset \M$ with ${\rm Volume} \, (\M \sim \M_\circ) < \e$}.
\end{align*}
In addition, if the local connection $1$-forms $\Omega_\U$ can be chosen to have higher regularity, then the  global connection $1$-form $\Omega$ also lies in the same regularity class.
\end{proposition}
Note that in Proposition~\ref{propn: luzin} we have assumed the local existence of connection given the curvature. The existence issue may be tricky in general: it amounts to a prescribing curvature problem. Nonetheless, in both of Uhlenbeck's classical works \cite{u', u} and our tentative applications to isometric immersions, we always start with a given, pre-existing family of connection $1$-forms $\{\Omega^\delta\}$ satisfying the equation $d\Omega^\delta+[\Omega^\delta \wedge \Omega^\delta] = \curv^\delta$, the goal is always to deduce properties of $\{\Omega^\delta\}$ (\emph{e.g.}, strong/weak compactness and regularity modulo gauge transforms) under certain assumptions on the associated curvature 2-forms $\curv^\delta$ (\emph{e.g.}, flat, Yang--Mills, or bounded in some $L^p$-norms, etc.).

\begin{remark}\label{remark, ppl G bdl}
All the results up to now carry over to any principal $G$-bundle $P$ with structure group $G$ over a closed Riemannian manifold $\M$. We may simply replace each occurrence of $\eta$ with $P$ and ${\rm Ad} \,\eta$ with $\lie$, respectively. Moreover, throughout the above we assume that $\M$ is compact, but it should be expected that variants of these results hold for open (noncompact $\&$ complete) Riemannian manifolds. In this case, one needs to invoke theories of elliptic PDE systems and Hodge decomposition on unbounded domains. This is left for future investigation.
 \end{remark}

\begin{remark}\label{rem: dropping []}
From now on, we shall write $$d\Omega + [\Omega \wedge \Omega] \equiv d\Omega + \Omega \wedge \Omega.$$ Similarly, we drop the  bracket for any other wedge products of $\lie$ (matrix Lie algebra)-valued differential 1-forms. This is a natural notational convention because for any $\lie$-valued 1-forms $A, B \in \Omega^1(\lie) = \G\left(\M;T^*\M \otimes \lie\right)$ one has the identity $$\left(A \wedge B\right)(u,v)=A(u)B(v)-B(u)A(v) \equiv [A \wedge B](u,v)$$ for all vector fields $u,v \in \G(T\M)$. This middle term is indeed an intertwining of Lie bracket (taken with respect to the $\lie$-factor) and wedge product of differential forms (taken with respect to the $T^*\M$-factor). 
\end{remark}

\subsection{Isometric immersions with low regularity}
The problem of isometric immersions and/or embeddings has been a crucial topic in global differential geometry for more than a hundred years, leading to seminal developments in geometric analysis --- \emph{e.g.}, Nash's $C^1$ and $C^{k}$-isometric embeddings ($k \geq 3$), Nash--Moser iteration, $h$-principle, Aleksandroff geometry, convex integration, and the analysis of fully nonlinear Monge--Amp\`{e}re equations, to name a few. For a detailed account of the history, status-quo, and open problems of isometric immersions, we refer to the comprehensive exposition \cite{hh} by Han--Hong and the problem section \cite{yau} by Yau.

As remarked at the very beginning, our note is essentially inspired by S. Mardare's works \cite{m05, m07}. These works are, in turn, largely motivated by the problem of the existence of isometric immersions of surfaces into Euclidean 3-space with $L^p_\loc$-second fundamental forms, under the assumption that the compatibility equations for curvatures --- \emph{i.e.}, the Gauss--Codazzi equations --- are satisfied in the distributional sense. This has been generalised to  isometric immersions of  arbitrary dimensions and codimensions with $L^p$-extrinsic geometry by Szopos \cite{sz}. Later, the author and Chen \cite{cl} gave a simplified geometric proof via the Cartan formalism. Recall too that the analogous result in the $C^\infty$-category is classical, under the name  \emph{fundamental theorem of surface theory}. See, \emph{e.g.}, Tenenblat \cite{ten}. An earlier extension with $C^0$-connections and second fundamental forms was obtained by Hartman--Wintner \cite{hw}. Moreover, the notable recent work \cite{l} by Litzinger extended the above results to the critical case $p=n=2$.\footnote{During the writing and revision of this work, the author and X. Su obtained a ``supercritical'' version of the fundamental theorem of surface theory in arbitrary dimensions and codimensions in \cite{lisu}. More precisely, for each $n \geq 3$, $W^{2,n}$-isometric immersions can be constructed over simply-connected $n$-dimensional domains from the weak solutions to the Gauss--Codazzi--Ricci equations in regularity class $L^n$. In fact, the $W^{2,n}$ and $L^n$ above can be replace with larges functions of the Morrey type, namely $L^{q,n-q}_2$ and $L^{q,n-q}$, respectively, for any $q \in ]2,n]$. }

A map  with Sobolev (\emph{e.g.}, $W^{2,p}$-) regularity, $$\iota: (\M,g)\map \R^{n+k},$$ is said to be an \emph{isometric immersion} iff $d\iota$ is one-to-one outside a null set of $\M$, and that $g$ coincides almost everywhere with the pullback of the Euclidean metric on $\R^{n+k}$ under $\iota$, namely that $d\iota \otimes d\iota = g$.  Throughout, we view $\R^{n+k}$ as equipped with the Euclidean metric. These aforementioned results obtained in \cite{m05, sz, cl} are summarised below:

\begin{proposition}\label{prop: isom imm}
Let $(\M,g)$ be a connected and simply-connected $n$-dimensional Riemannian manifold with $W^{1,p}_\loc$-metric; $p>n$. Suppose that the second fundamental form $\two$ and the normal affine connection $\na^{\perp}$ (\emph{i.e.}, the orthogonal projection of the Levi-Civita connection on $\R^{n+k}$ on the normal bundle) have $L^p_\loc$-components in local coordinates, and that they together with $g$ satisfy the Gauss--Codazzi--Ricci equations in the sense of distributions. Then there exists a locally $W^{2,p}$-isometric immersion $\iota: (\U\subset\M,g) \emb \R^{n+k}$ such that $\left(\two, \na^{\perp}\right)$ coincides with its extrinsic geometry. In addition, $\iota$ is unique modulo Euclidean rigid motions in  $\R^{n+k}$ up to null sets.
\end{proposition}


Proposition~\ref{prop: isom imm} above is heuristic, just to give the key ideas. One issue was left ambiguous: without \emph{a priori} knowledge on the existence of isometric immersion, what does it mean by the second fundamental form $\two$ and the normal affine connection $\na^\perp$, not to mention that they satisfy (weakly) the Gauss--Codazzi--Ricci equations?

To address this issue, rigorous versions of Proposition~\ref{prop: isom imm} are presented below as Theorems~\ref{thm: existence, early} and \ref{thm: existence}. The idea is: besides the manifold $(\M,g)$, our data includes a given abstract vector bundle $\pi: E \to \M$ of rank $k$ with some Riemannian metric $g^E$. The second fundamental form $\two$ associated to $E$ is formally given by the tensorial rules one should have if $E$ is indeed the normal bundle of an isometric immersion, and the  connection $\na^E$ is metric-compatible with $g^E$. Then, one may also write the Gauss, Codazzi, and Ricci equations associated to $E$ by splitting the zero Riemann curvature tensor into $(T\M, T\M)$, $(T\M, E)$, and $(E,E)$-factors, respectively. The conclusion of Proposition~\ref{prop: isom imm}, if rigorously interpreted, consists of that the normal bundle of $\iota$ can be identified with $E$ via a bundle isomorphism.

Before presenting the theorem, let us fix once and for all our convention for indices:
\begin{equation}\label{indices, convention}
1\leq i,j,k\leq n;\quad n+1\leq\alpha,\beta \leq n+k; \quad 1 \leq a,b,c \leq n+k.
\end{equation}
Here $i,j,k$ are indices for the tangent bundle $T\M$, and $\alpha,\beta$ are indices for the bundle $E$. Accordingly, we shall also use capital Latin (Greek, resp.) letters to denote vectorfields on $\M$ (sections of $E$, resp.). For a vector bundle $E$ we write $\G(E)$ for the space of its sections, whose regularity assumption is clear from the context.

\begin{theorem}\label{thm: existence, early}
Let $(\M,g)$ be an $n$-dimensional Riemannian manifold with metric $g \in W^{1,p}_\loc$; $p>n$. Let $E$ be a vector bundle of rank $k$ over $\M$. Assume $E$ is equipped with a $W^{1,p}_\loc$-metric $g^E$ and an $L^p_\loc$-connection $\na^E$ compatible with $g^E$. Given an $L^p$-tensor field
\begin{equation*}
\mathcal{S}: \G(E) \times \G(T\M) \longrightarrow \G(T\M), \qquad \mathcal{S}(\eta, X) \equiv  \mathcal{S}_\eta X 
\end{equation*}
such that 
\begin{equation}
g\big(X, \mathcal{S}_\eta Y\big) = g\big( \mathcal{S}_\eta X, Y\big)
\end{equation}
for all $X,Y \in \G(T\M)$ and $\eta \in \G(E)$. Then, define $$\two:\G(T\M) \times \G(T\M) \map \G(E)$$ via
\begin{equation}
g^E\big(\two(X,Y),\eta\big) := -g\big(\mathcal{S}_\eta X, Y\big).
\end{equation}

The following are equivalent:
\begin{enumerate}
\item[(I)]
For each $x \in \M$, there exists an open neighbourhood $\U \subset \M$ containing $x$ and  a \emph{local} isometric immersion $\iota: (\U,g) \map \R^{n+k}$ of $W^{2,p}$-regularity, whose normal bundle $T\R^{n+k}\slash T(\iota \M)$, Levi-Civita connnection on the normal bundle, and second fundamental form can be identified with $E$, $\na^E$, and $\two$, respectively. In addition, $\iota$ is  unique up to the Euclidean rigid motions in $\R^{n+k}$ modulo null sets.
\item[(II)]
The \emph{Cartan formalism} holds in the sense of distributions:
\begin{eqnarray}
&&d\omega^i = \sum_j \omega^j \wedge \Omega^i_j;\label{first structure eq}\\
&&0=d\Omega^a_b + \sum_c \Omega^c_b\wedge \Omega^a_c,\label{second structure eq}
\end{eqnarray}
where $\{\omega^i\}_{1\leq i \leq n}$ is an orthonormal coframe for $(T^*\M,g)$, and $\{\Omega^a_b\}_{1\leq a,b \leq n+k}$ is the connection $1$-form given by
\begin{eqnarray}
&& \Omega^i_j (\p_k) := g(\na_{\p_k}\p_i,\p_j);\label{Omega 1}\\
&& \Omega^i_\alpha (\p_j) \equiv -\Omega^\alpha_i (\p_j) := g^E\big(\two(\p_i,\p_j), \eta_\alpha\big);\label{Omega 2}\\
&& \Omega^\alpha_\beta(\p_j):=g^E\big(\na^E_{\p_j}\eta_\alpha,\eta_\beta\big).\label{Omega 3}
\end{eqnarray}
In the above, $\{\p_i\}$ is the orthonormal frame for $(T\M,g)$ dual to $\{\omega^i\}$, and $\{\eta_\alpha\}_{n+1\leq \alpha \leq n+k}$ is an orthonormal frame for $(E,g^E)$.
\item[(III)]
The Gauss--Codazzi--Ricci equations hold in the sense of distributions:
\begin{eqnarray}
&& g\big(\two(X,Z), \two (Y,W)\big) - g\big( \two(X,W),\two(Y,Z) \big) = R(X,Y,Z,W);\label{gauss}\\
&& \overline{\na}_Y\two(X,Z) - \overline{\na}_X\two(Y,Z)=0;\label{codazzi}\\
&& g\big([\mathcal{S}_\eta, \mathcal{S}_{\zeta}] X, Y\big) = R^E(X,Y,\eta,\zeta)\label{ricci},
\end{eqnarray}
for all $X,Y,Z,W \in \G(T\M)$ and $\eta,\zeta \in \G(E)$. Here, $[\bullet,\bullet]$ is the commutator of operators,  $R$ and $R^E$ are respectively the Riemann curvature tensors for $(T\M,g)$ and $\left(E,g^E\right)$, and $\overline{\na}$ is the Levi-Civita connection on Euclidean space $\R^{n+k}$.
\end{enumerate}
\end{theorem}

From the PDE perspectives, in Eqs.~\eqref{gauss}, \eqref{codazzi}, and \eqref{ricci} the metric $g$ (\emph{i.e., the intrinsic geometry}; hence the Levi-Civita connection $\na$ and the Riemann curvature tensor $R$) is given and $\two$, $\na^E$ (\emph{i.e., the extrinsic geometry}) are unknown.
Eqs.~\eqref{first structure eq} and \eqref{second structure eq} are known as the \emph{first} and the \emph{second structural equations} of Cartan's formalism for Euclidean submanifolds, respectively. Note that the second  structural equation is nothing but Eq.~\eqref{curvature, def} with trivial curvature 2-form.

The structural equations~\eqref{first structure eq} and \eqref{second structure eq}, written component-wise in local co-ordinates as such, are identities in $\Omega^2(M)=\G\left(T^*\M^{\otimes 2}\right)$ for fixed indices $i,a,b$. Thus, these equations should be understood via the following identities:   $\left(A \wedge B\right)(u,v)=A(u)B(v)-B(u)A(v)$ (recall Remark~\ref{rem: dropping []}) and $dA(u,v)=u(A(v))-v(A(u))-A([u,v])$ for vector fields $u,v$ and 1-forms $A,B$.

One natural question arises: when can the $W^{2,p}$-local isometric immersion $\iota$ in Theorem~\ref{thm: existence, early} (I) be made \emph{global}? A classical sufficient condition is that $\M$ is simply-connected (\emph{cf}. Tenenblat \cite[Remark on p.35]{ten}; see also Ciarlet--Larsonneur \cite{ciarlet'}). We give another sufficient condition below: \emph{$L^p$-smallness of $\Omega$}. This is equivalent to the $L^p$-smallness of both the Levi-Civita connection of $(\M,g)$ and the extrinsic geometry, \emph{i.e.}, the second fundamental form $\two$ and the connection $\na^E$ on the vector bundle $E$. 

\begin{theorem}\label{thm: existence}
In the setting of Theorem~\ref{thm: existence, early}, let $\M$ be a complete Riemannian manifold and assume (II) or (III) --- \emph{i.e.}, assume either the Cartan formalism or the Gauss--Codazzi--Ricci equations hold in the sense of distributions. Then, under either of the following additional assumptions, the local isometric immersion $\iota$ in (I) can be promoted to a \emph{global} isometric immersion $\iota: (\M,g) \to \R^{n+k}$ of $W^{2,p}_\loc$-regularity:
\begin{enumerate}
\item
$(\M,g)$ is simply-connected (compact or noncompact).
\item
$(\M,g)$ is a closed manifold which is ``nearly flat'', in the sense that the connection 1-form $\Omega$ (defined via Eqs.~\eqref{Omega 1}--\eqref{Omega 3} in Theorem~\ref{thm: existence, early} (II)) lies in $L^p\left(\M; T^*\M\otimes \so(n+k,\R)\right)$, with $\|\Omega\|_{L^p(\M)} <$ the uniform constant $\kappa_0(p,\M)$ as in Theorem~\ref{thm: main}.

\end{enumerate}
\end{theorem}




We shall present  in \S\ref{sec: isom imm} a  proof of Theorem~\ref{thm: existence}, which is essentially a recast of the arguments in \cite[\S 5]{cl}. We utilise ideas from Tenenblat \cite{ten} and S. Mardare \cite{m05}, as well as applying  Theorem~\ref{thm: main} proved earlier in this note. Our intention here is to elucidate some issues on local \emph{vs.} global isometric immersions in the literature. Note also that Theorem~\ref{thm: existence} (2) is a new result elusive in the existing literature. 

Without either the simple connectedness of $\M$ or the $L^p$-smallness condition on $\Omega$, at the moment we do not know whether global $W^{2,p}_\loc$-isometric immersions always exist given the Cartan formalism/Gauss--Codazzi-Ricci equations. Nevertheless, an ``almost global''  $W^{2,p}_\loc$-isometric immersions can be easily established. ``Almost global'' is again understood in the Luzin sense.

\begin{theorem}\label{thm: luzin, almost global isom imm}
In the setting of Theorem~\ref{thm: existence, early}, let $\M$ be a closed Riemannian manifold and assume (II) or (III) --- assume either the Cartan formalism or the Gauss--Codazzi--Ricci equations hold in the sense of distributions. Then for any $\e>0$ there is an closed set $\M_\times\subset\M$ such that 
\begin{itemize}
\item
${\rm Volume}\left(\M_\times\right)<\e$; and 
\item
there exists a $W^{2,p}$-isometric immersion $\iota:\left(\M\sim\M_\times,g\right)\emb \R^{n+k}$.
\end{itemize}
\end{theorem}

Discussion will also be made concerning restrictions imposed on the connection 1-forms for metrics \emph{without} smooth isometric immersions (see Pogorelov \cite{p}, Iaia \cite{iaia}, and Khuri \cite{khuri}, among other works) in the framework of the Cartan formalism of structural equations, from the perspectives of smoothability of connections as in Theorem~\ref{thm: main}.

 Towards the end of \S\ref{sec: isom imm}, we briefly discuss the weak continuity of the Gauss--Codazzi--Ricci equations. Whenever $p>2$, \emph{regardless of dimension $n$ and codimension $k$}, weak  solutions to the Cartan structural equations (or equivalently, to the Gauss--Codazzi--Ricci equations) are well defined, and the weak continuity property of these equations has been established by Chen--Slemrod--Wang \cite{csw} (see also \cite{cl, cg, giron}). We shall point out that essentially the same arguments carry over to isometric immersions into general ambient manifolds.

\subsection{Organisation}

The remaining parts of the paper is organised as follows.

In \S\ref{sec: smoothability} we shall prove our main Theorem~\ref{thm: main} on smoothability of $L^p$-connections with prescribed curvature. In \S\ref{sec: existence} we prove Corollary~\ref{cor: regularity} and Proposition~\ref{propn: luzin}, which are direct consequences of  Theorem~\ref{thm: main}. Next, in \S\ref{sec: isom imm} we discuss the application of our smoothability  Theorem~\ref{thm: main} to $W^{2,p}$-isometric immersions. Finally, we conclude the note with remarks in \S\ref{sec: remarks}.

\section{Smoothability of $L^p$-connections} \label{sec: smoothability}

This section is devoted to the proof of Theorem~\ref{thm: main}, which is rephrased as follows. Throughout the remaining parts of the paper, we  adhere to the notations in Remark~\ref{rem: dropping []}.

\begin{theorem}\label{thm: re}
Let $\eta$ be a vector bundle with structure group $G$ over a closed Riemannian $n$-dimensional manifold $\M$; $n \geq 2$. 
Assume that the connection 1-form $\Omega$ lies in $\conn^{p}_0$ for some $p>n$ and that the associated curvature 2-form $\curv:=d\Omega + \Omega \wedge \Omega$ is in the regularity class $W^{s,q}$, where $s \geq 0$ and $\frac{p}{2} \leq q \leq \infty$ such that $W^{s,q}\subset L^2$ over $\M$. Then there is a uniform constant $\kappa_0>0$ depending only on $p$ and the geometry of $\M$, such that if
$$\|\Omega\|_{L^{p}(\M)} < \kappa_0,$$ 
then one can find a one-parameter family $\{\Omega^\e\}_{0<\e<1} \subset \conn_{1+s}^q$ such that $\Omega^\e \to \Omega$ strongly in  $L^p$, and that the curvature of $\Omega^\e$ coincides with $\curv$ for each $\e$ on $\M$.
\end{theorem}



We also have the following local version of the above smoothability result: 

\begin{theorem}\label{thm: local}

Let $\eta$, $G$ be as in Theorem~\ref{thm: re} and $\M$ be an $n$-dimensional manifold (not necessarily compact). Suppose that the connection 1-form $\Omega$ and its associated curvature 2-form $\curv$ are in $L^p_\loc$ and $W^{s,q}_\loc$ over $\M$, respectively, with $p,s,q$ as in Theorem~\ref{thm: re}. Then, there exists a uniform constant $\kappa_0>0$ depending only on $p$ and $\M$, such that the following holds. 

For any bounded smooth domain $U \Subset \M$, if $\|\Omega\|_{L^p(U)}<\kappa_0$, then one can find a one-parameter family $\{\Omega^\e\}_{0<\e<1}$ of connection 1-forms on $U$ with $W^{1+s,q}$-regularity up to the boundary, such that $\Omega^\e \to \Omega$ in $L^p(U)$ and $F_{\Omega^\e}\big|_U := \left(d\Omega^\e+\Omega^\e \wedge \Omega^\e\right)\big|_U = \curv\big|_U$ for each $\e$. 
\end{theorem}

In the remaining part of this section we shall prove Theorems~\ref{thm: re} and \ref{thm: local} at one stroke. The proof of Theorem~\ref{thm: re} essentially follows from (and is somewhat simpler than) that of its local variant, \emph{i.e.}, Theorem~\ref{thm: local}, as we do not need to tackle with boundary conditions over closed manifolds. Nevertheless, in this case we need to pay extra care to fix a cohomology class. 

Throughout, $\ttt$ and $\nnn$ denote respectively the tangential and normal trace operators. For a Riemannian manifold-with-boundary $\M$, write $\iota: \p\M\emb\M$ for the inclusion of the boundary $\p\M$ into $\M$. Then, for any differential form $\alpha \in \Omega^\bullet(\M)$, set 
\begin{equation*}
\ttt\alpha := \iota^\#\alpha,\qquad \nnn\alpha:=\alpha - \ttt\alpha,
\end{equation*}
both being differential forms defined over $\p\M$. See Schwarz \cite[Introduction, p.2]{sch} for details.

\begin{proof}[Proof of Theorems~\ref{thm: re} and \ref{thm: local}]
The arguments are divided into nine steps. In Steps~1--8 we prove Theorem \ref{thm: local} for the special case $\M=\R^n$. Then, in Step~9 we extend to the case of general closed Riemannian manifolds and conclude Theorems~\ref{thm: re} and \ref{thm: local}.

If $\M=\R^n$, one may globally identify $\bigwedge^2 T^*\R^n$ with $\so(n,\R)$, the space of $n \times n$ antisymmetric real matrices. For $k \in \mathbb{N}$, a harmonic $k$-form $\harm$ satisfies that $d\harm = 0$ and $d^*\harm=0$.

\smallskip
\noindent
{\bf Step 1.} We first establish:

\noindent
{\bf Claim A.} Let $U \Subset \R^n$ be a smooth Euclidean domain. There are a scalar field $\phi \in W^{1,p}\left(U; {\rm Ad}\,\eta\right)$, a differential 2-form $\psi \in W^{1,p}\left(U; {\rm Ad}\,\eta  \otimes \bigwedge^2 T^*\R^n\right)$ with $\nnn\psi = 0$ and $\nnn d\psi=0$ on $\p U$, and a harmonic 1-form $\harm \in C^\infty\left(U; {\rm Ad}\,\eta  \otimes \bigwedge^2 T^*\R^n\right)$ with $\nnn \harm = 0$ on $\p U$, satisfying that 
\begin{equation}\label{new, decomposition of Omega}
\Omega = d\phi + d^*\psi +\harm \qquad \text{ in } U.
\end{equation}
The three terms on the right-hand side of Eq.~\eqref{new, decomposition of Omega} are $L^2$-orthogonal to each other.

\begin{proof}[Proof of Claim A]

This follows from classical Hodge decomposition on manifolds-with-boundary. We first solve for $\psi$ from the boundary value problem:
\begin{equation}\label{psi system}
\begin{cases}
\Delta \psi = \curv - \Omega \wedge \Omega\qquad \text{ in } U,\\
\nnn\psi = 0 \quad \text{ and } \quad \nnn d\psi=0 \quad \text{ on } \p U,
\end{cases}
\end{equation}
where $$\Delta:=dd^*+d^*d$$ is the Laplace--Beltrami operator. For $p>n$, by  Sobolev embedding  $W^{1,p'}_\loc \emb L^{\frac{p}{{\frac{n-1}{n}\cdot p}-1}}_\loc \emb L^{\frac{p}{p-2}}_\loc$, and hence by duality  $L^{\frac{p}{2}}_\loc=L^{\left(\frac{p}{p-2}\right)'}_\loc \emb W^{-1,p}_\loc$. All the Sobolev spaces here are defined over $\R^n$. Then we can solve Eq.~\eqref{psi system} --- which is elliptic in the sense of Agmon--Douglis--Nirenberg \cite{adn1, adn2} or \u{S}apiro--Lopatinski\u{i} --- from the standard elliptic theory to get a solution $\psi$ in $W^{1,p}$ (\cite{sch}). 

Moreover, since $d(\curv-\Omega \wedge \Omega)=dd\Omega=0$ and $d\Delta=\Delta d$, we have 
\begin{align*}
\Delta d\psi=0\qquad \text{ in } U.
\end{align*}
This together with the boundary condition in Eq.~\eqref{psi system} shows that $\psi$ is closed. Then 
\begin{align*}
d(\Omega-d^*\psi) = d\Omega-\Delta\psi = d\Omega + \Omega \wedge \Omega - \curv = 0\quad \text{ in } U,
\end{align*}
where the penultimate equality follows from the definition of $\psi$ via Eq.~\eqref{psi system}.

To conclude the proof, we apply the version of Friedrichs--Hodge--Morrey decomposition as in \cite[Corollary~2.4.9]{sch} to decompose 
\begin{equation}\label{new, friedrichs}
\Omega - d^*\psi = d\alpha_0 + d^*\beta_0 + d\gamma_0 + \lambda_0,
\end{equation}
where $\alpha_0$, $\beta_0$ are of regularity $W^{1,p}$, and $d\gamma_0$, $\lambda_0$ are harmonic 1-forms where $\nnn \lambda_0 = 0$ on $\p U$. In addition, the four components on the right-hand side are $L^2$-orthogonal with respect to the usual $L^2$-inner product on differential forms. That is,
\begin{align*}
\bra\alpha,\beta\ket_U := \int_U \alpha \wedge \star\beta
\end{align*}
with $\star$ being the Hodge star operator associated to the Riemannian volume measure $U \subset \M$. The wedge product $\wedge$ is intertwined with the product on the Lie algebra $\lie$, if $\alpha$ and $\beta$ are $\lie$-valued differential forms (here $\lie = \so(n;\R)$). Since $d(\Omega - d^*\psi)=0$, we test Eq.~\eqref{new, friedrichs} against $\beta_0$ and integrate over $U$ to obtain via orthogonality that
\begin{align*}
0 &= \bra d\left(\Omega - d^*\psi\right),\, \beta_0\ket \\
&= \bra \Omega - d^*\psi,\, d^*\beta_0\ket \\
&= \bra d\alpha_0 + d^*\beta_0 + d\gamma_0 + \lambda_0,\, d^*\beta_0\ket\\
&= \bra d^*\beta_0,\, d^*\beta_0\ket.
\end{align*}
Thus $d^*\beta_0=0$. We may rewrite Eq.~\eqref{new, friedrichs} as
\begin{equation*}
\Omega - d^*\psi = d\phi + \harm
\qquad \text{ in } U,
\end{equation*}
where $\psi = \alpha_0+\gamma_0$ and $\lambda_0=\harm$. Hence \emph{Claim A} follows.   \end{proof}

\smallskip
\noindent
{\bf Step 2.} Next, we take $\{\phi^\e\} \subset C^\infty(U; {\rm Ad}\,\eta )$ such that
\begin{equation}\label{phi epsilon}
\phi^\e \longrightarrow \phi\qquad \text{ in } W^{1,p}\left(U;{\rm Ad}\,\eta \right).
\end{equation} 
When for $p>n$, such $\{\phi^\e\}$ always exists because $\phi$ has a H\"{o}lder continuous representative. By isometrically embedding ${\rm Ad}\,\eta$ into a Euclidean space via Nash's theorem as a submanifold, one may  construct $\phi^\e$ via standard mollifications followed by reprojections onto the image.

Note that we can take here $\{\phi^\e\} \subset C^\infty\left(\overline{U}; {\rm Ad}\,\eta \right)$, \emph{i.e.}, smooth up to the boundary. See, for instance, Evans \cite[p.252, 5.3.3, Theorem~3]{evans}.

\smallskip
\noindent
{\bf Step 3.} To proceed, let us solve for $\psi^\e$, smooth approximates of $\psi$,  from the following nonlinear elliptic system ($\harm$ is the harmonic 1-form in Step~1, Eq.~\eqref{new, decomposition of Omega}):
\begin{equation}\label{system for psi epsilon}
\begin{cases}
\Delta\psi^\e = \curv - \left(d\phi^\e + d^* \psi^\e + \harm\right) \wedge \left(d\phi^\e + d^* \psi^\e+\harm\right) \qquad \text{ in } U,\\
\nnn\psi^\e = 0 \quad \text{ and } \quad \nnn d\psi^\e=0 \quad \text{ on } \p U.
\end{cases}
\end{equation}
In this and the next steps, we shall solve Eq.~\eqref{system for psi epsilon} via Schauder's fixed point argument.

For this purpose, denote by $\Delta^{-1}$ the solution operator for  Laplace--Beltrami subject to the same boundary conditions. Consider the operator for each fixed $\e\in]0,1]$:
\begin{align}\label{operator T}
\T_\e:\quad W^{1,p}\left(U;  {\rm Ad}\,\eta \otimes \bigwedge^2 T^*\R^n\right) &\longrightarrow W^{1,p}\left(U; {\rm Ad}\,\eta \otimes \bigwedge^2 T^*\R^n\right)\nonumber\\
 \zeta \quad &\longmapsto\quad \Delta^{-1}\Big\{ \curv-\left(d\phi^\e + d^*\zeta+\harm\right)\wedge\left(d\phi^\e + d^*\zeta+\harm\right) \Big\}.
\end{align}

\noindent
{\bf Claim B.} For each $0<\e\leq 1$, the operator $\T_\e$ in Eq.~\eqref{operator T} is compact.

\begin{proof}[Proof of Claim B] We first observe that $\T_\e$ indeed maps into $W^{1,p}\left(U;  {\rm Ad}\,\eta \otimes \bigwedge^2 T^*\R^n\right)$. This is because for $\zeta \in W^{1,p}\left(U;{\rm Ad}\,\eta  \otimes \bigwedge^2 T^*\R^n\right)$, by Sobolev embedding one has  $$\curv-\left(d\phi^\e + d^*\zeta + \harm\right)\wedge\left(d\phi^\e + d^*\zeta +\harm \right) \in W^{-1,p}\left(U;{\rm Ad}\,\eta \otimes \bigwedge^2 T^*\R^n\right).$$ Hence, by standard elliptic estimates, $\T_\e\zeta \in W^{1,p}\left(U;{\rm Ad}\,\eta  \otimes \bigwedge^2 T^*\R^n\right)$. 

In addition, we can easily bound
\begin{align*}
\|\T_\e\zeta\|_{W^{1,p}(U)} &\leq C\left\|\curv-\left(d\phi^\e + d^*\zeta+\harm\right)\wedge\left(d\phi^\e + d^*\zeta+\harm\right) \right\|_{W^{-1,p}(U)}\\
&\leq C \left(\left\|\curv\right\|_{W^{-1,p}(U)} +\left\|\left(d\phi^\e + d^*\zeta+\harm\right)\wedge\left(d\phi^\e + d^*\zeta+\harm\right) \right\|_{L^{\frac{np}{n+p}}(U)} \right)\\
&\leq C \left(\left\|\curv\right\|_{W^{-1,p}(U)} +\left\|d\phi^\e + d^*\zeta +\harm\right\|^2_{L^{\frac{2np}{n+p}}(U)} \right)\\
&\leq C\left(\left\|\curv\right\|_{W^{-1,p}(U)} +\left\|d\phi^\e + d^*\zeta+\harm \right\|^2_{L^{p}(U)} \right)\\
&\leq C\left(\left\|\curv\right\|_{W^{-1,p}(U)} +\left\|\phi^\e\right\|^2_{W^{1,p}(U)}+\left\|\zeta\right\|^2_{W^{1,p}(U)} + \|\harm\|^2_{L^p(U)} \right),
\end{align*}
where $C$ depends only on $p$. The first line follows from standard elliptic estimates, the second line holds by the Sobolev embedding $W^{1,\frac{np}{n+p}}\emb L^p$ for $p>n$ over $U \Subset \R^n$, the third line follows from the Cauchy--Schwarz inequality, and the fourth line holds by $p\geq n$ and the boundedness of $U$. Thus, in view of the definition of $\psi^\e$ in Eq.~\eqref{phi epsilon}, we have 
\begin{equation}\label{W1,p estimate for solution operator}
\|\T_\e\zeta\|_{W^{1,p}(U)} \leq C(p) \left(\left\|\curv\right\|_{W^{-1,p}(U)} +\left\|\phi\right\|^2_{W^{1,p}(U)}+\left\|\zeta\right\|^2_{W^{1,p}(U)} + \|\harm\|_{L^p(U)}^2+\mathfrak{o}(\e) \right)
\end{equation}
with $\mathfrak{o}(\e) \to 0$ as $\e \to 0$. In addition, recall from the assumption that $\curv \in W^{s,q}\left(\M;\lie \otimes \bigwedge^2 T^*\M\right)$ with $s \geq 0$, $q \geq \frac{p}{2}$, as well as the compact Sobolev embedding $W^{s,q}\emb\emb W^{-1,p}$. Hence, Eq.~\eqref{W1,p estimate for solution operator} yields that  $\T_\e$ is $W^{1,p} \to W^{1,p}$ bounded uniformly in $\e$.


Next we check that $\T_\e$ is a compact operator. Take any sequence $\left\{\zeta^\e_j\right\}_{j=1,2,\ldots}$ of $\so(n;\R)$-valued 2-forms on $U$ with $W^{1,p}$-regularity, which converges weakly in $W^{1,p}$ to $\overline{\zeta}^\e$. Then
\begin{align*}
&\Big\{\curv-\left(d\phi^\e + d^*\zeta_j^\e+\harm\right)\wedge\left(d\phi^\e + d^*\zeta^\e_j+\harm\right) \Big\}\\
&\qquad\qquad - \Big\{\curv-\left(d\phi^\e + d^*\overline{\zeta}^\e + \harm\right)\wedge\left(d\phi^\e + d^*\overline{\zeta}^\e+\harm\right) \Big\}\\
&\quad = d^*\left(\overline{\zeta}^\e-\zeta_j^\e\right) \wedge \left(d\phi^\e + d^*\overline{\zeta}+\harm\right) +  \left(d\phi^\e + d^*{\zeta}_j^\e + \harm\right) \wedge d^*\left(\overline{\zeta}^\e-\zeta_j^\e\right)
\end{align*}
is uniformly bounded in $L^{p/2}\left(U;{\rm Ad}\,\eta \otimes \bigwedge^2 T^*\R^n\right)$ by Cauchy--Schwarz and that $\phi^\e \to \phi$ in $W^{1,p}$ (Step~2, Eq.~\eqref{phi epsilon}). Note here that the second term on the right-hand side is a pairing of weakly $L^p$-convergent sequences.

To proceed, under the assumption $p>n$ we have the compact Sobolev embedding:
\begin{align}\label{compact embedding}
 L^{\frac{p}{2}}\left(U;{\rm Ad}\,\eta\otimes\bigwedge^2 T^*\R^n\right) \emb \emb  W^{-1,{p}}\left(U;{\rm Ad}\,\eta\otimes\bigwedge^2 T^*\R^n\right)
\end{align}
Thus, one may select a subsequence $\left\{\zeta^\e_{j_\nu}\right\} \subset \left\{\zeta^\e_j\right\}$ such that $\zeta_{j_\nu}^\e \to \overline{\zeta}^\e$ strongly in the $W^{-1,p}$-topology as $\nu \to \infty$. By standard elliptic estimates applied to Eq.~\eqref{system for psi epsilon}, we have
\begin{align*}
\lim_{\nu \to \infty}\left\|\T_\e \zeta_{j_\nu}^\e - \T_\e\overline{\zeta}^\e\right\|_{W^{1,p}(U)} =  0. 
\end{align*}
Thus the operator $\T_\e$ defined in Eq.~\eqref{operator T} is compact.

As a consequence, by Schauder's fixed point theorem, for each $0<\e\leq 1$ there exists a solution $\psi^\e \in W^{1,p}\left(U; {\rm Ad}\,\eta  \otimes \bigwedge^2 T^*\R^n\right)$ for Eq.~\eqref{psi system}. Hence \emph{Claim B} is proved.  \end{proof}

Note that the solution operator $\T_\e$ in Eq.~\eqref{operator T} is indeed uniformly $W^{1,p}\to W^{1,p}$ bounded in $\e$. That is, there is a uniform constant $K_0$ independent of $\e$ such that
\begin{equation}\label{unif bdd of solution operator}
\sup_{0< \e \leq 1}\||\T_\e\|| \leq K_0,
\end{equation}
where $\||\bullet\||$ is the strong $W^{1,p}\to W^{1,p}$ operator norm. This follows directly from Eq.~\eqref{W1,p estimate for solution operator}. 

\smallskip
\noindent
{\bf Step 4.} Now we establish that the operator $\T_\e$ indeed has a unique fixed point, under the $L^p$-smallness condition~\eqref{smallness condition--key} on $\Omega$. (We are unable to show the uniqueness of the solution $\psi^\e$ to Eq.~\eqref{psi system} for a general $L^p$-connection $\Omega$, even when the connection is smooth.) 

\smallskip
\noindent
{\bf Claim C.} For each $0<\e\leq 1$, the operator $\T_\e$ has a fixed point $\psi^\e$, provided that 
\begin{equation}\label{smallness condition--key}
\|\curv\|_{W^{-1,p}(U)} + \|\Omega\|_{L^p(U)} \leq \kappa
\end{equation}
for a small uniform number $\kappa$ depending only on $p$. 

Before presenting the proof, let us remark that the smallness of $\|\curv\|_{W^{-1,p}(U)}$ follows from that of $\|\Omega\|_{L^p(U)}$, so the above hypothesis~\eqref{smallness condition--key} is equivalent to the condition in the statement of Theorem~\ref{thm: main}. 
Indeed, $\|\Omega\|_{L^p(U)}\leq \kappa$ implies that $$\|d\Omega\|_{W^{-1,p}(U)}\leq\kappa$$ and that
\begin{align}
\|\Omega\wedge\Omega\|_{W^{-1,p}(U)}\leq C' \|\Omega\wedge\Omega\|_{L^{\frac{p}{2}}(U)}\leq  C'\|\Omega\|_{L^p(U)}^2\ \leq C'\kappa^2.\label{C', new}
\end{align} Here $C'$ is the Sobolev constant corresponding to $$L^{\frac{p}{2}}\left(U;{\rm Ad}\,\eta\otimes\bigwedge^2 T^*\R^n\right) \emb\emb  W^{-1,p}\left(U;{\rm Ad}\,\eta\otimes\bigwedge^2 T^*\R^{n}\right),$$ which can be chosen to depend only on $p$ and $n$. The compact embedding holds whenever $p>n$. See the proof of \emph{Claim~A} in Step~1.

\begin{proof}[Proof of Claim~C]
Denote by $C_1=C(p)$ the uniform constant in Eq.~\eqref{W1,p estimate for solution operator} in the proof of \emph{Claim B} above, and by $C_2$ the constant $C'$ in Eq.~\eqref{C', new}.

We first fix a $\kappa_1>0$ so small that 
\begin{align*}
C_1 (\kappa_1)^2 \leq \frac{\kappa_1}{3}.
\end{align*}
In view of Eq.~\eqref{C', new} and the definition of curvature 2-form in Eq.~\eqref{curvature, def}, we have $\|\curv\|_{W^{-1,p}(U)} \leq \kappa + C_2 \kappa^2$. Meanwhile, by \emph{Claim~A}, Eq.~\eqref{new, decomposition of Omega} we have  $\|\harm\|_{L^p(U)}^2 + \|\phi\|^2_{W^{1,p}(U)} \leq \kappa^2$. By passing to a subsequence in Eq.~\eqref{phi epsilon} if necessary, we may also assume that $
\mathfrak{o}(\e)=\|\phi^\e - \phi\|_{W^{1,p}(U)} \leq \kappa$ for all $0<\e \leq 1$. Now, choose $\kappa$ depending only on $\kappa_1$, $C_1$, and $C_2$ (hence $\kappa= \kappa(p,n)$) such that
\begin{align*}
C_1\big[ \left(2+C_2\right)\kappa^2 + 2\kappa \big] \leq \frac{\kappa_1}{6}.
\end{align*}
Hence,
\begin{align*}
C_1\big[ \left(2+C_2\right)\kappa^2 + 2\kappa + (\kappa_1)^2 \big] \leq \frac{\kappa_1}{2}.
\end{align*}
By virtue of Eq.~\eqref{W1,p estimate for solution operator}, this implies that the restricted operator $\T_\e\big|_{
\mathfrak{C}_{\kappa_1}}: 
\mathfrak{C}_{\kappa_1} \to
\mathfrak{C}_{\kappa_1}$ is coercive under the smallness condition~\eqref{smallness condition--key}. That is, the $W^{1,p}\to W^{1,p}$-operator norm 
\begin{equation}\label{coercivity}
\sup_{0<\e\leq 1}\Big|\Big\|\T_\e\big|_{
\mathfrak{C}_{\kappa_1}}\Big\|\Big| <1
\end{equation}
whenever Eq.~\eqref{smallness condition--key} is valid. 

The existence and uniqueness of the fixed point $\psi^\e$ now follows from Eq.~\eqref{coercivity} and the Banach fixed point theorem. Hence \emph{Claim C} is proved.  \end{proof}

\smallskip
\noindent
{\bf Step 5.} In Steps~3 and 4 above, we have proved the unique solubility of Eq.~\eqref{system for psi epsilon} for $\psi^\e$ under the smallness assumption~\eqref{smallness condition--key}. The regularity of $\psi^\e$ follows directly from a bootstrap argument as presented below.

\noindent
{\bf Claim D.} Assume that the prescribed curvature $2$-form $\curv$ has $W^{s,q}$-regularity for $s \geq 0$ and $q \geq \frac{p}{2}>\frac{n}{2}$, such that $W^{s,q}\subset L^2$ over $U$. Then we have
\begin{align}\label{psi, regularity}
\psi^\e \in W^{s+2,q}\left(U; {\rm Ad}\,\eta \otimes \bigwedge^2T^*\R^n\right),
\end{align}
with $\|\psi^\e \|_{W^{s+2,q}(U)} \leq K_2 = K_2(n,p)$ uniformly in $0<\e \leq 1$. The regularity of $\psi^\e$ persists up to the boundary $\p U$.

\begin{proof}[Proof of Claim D]
Starting with a solution $\psi^\e$ for Eq.~\eqref{system for psi epsilon} in Step~3, we find that the term on the right-hand side, $\curv- \left(d\phi^\e + d^* \psi^\e+\harm\right) \wedge \left(d\phi^\e + d^*\psi^\e+\harm\right),$ is an $\so(n;\R)$-valued 2-form over $U\Subset \R^n$ of regularity $W^{s,q}+L^{p/2}$. By assumption, $s \geq 0$ and $q \geq p/2$, so it is of regularity $L^{p/2}$. Standard Elliptic regularity theory ($W^{2,p}$-estimates of Calder\'{o}n--Zygmund) applied to Eq.~\eqref{system for psi epsilon} yields $\psi^\e \in W^{2,{p}/{2}}\left(U; {\rm Ad}\,\eta \otimes \bigwedge^1T^*\R^n\right)$, which implies that
\begin{small}
\begin{align*}
\curv - \left(d\phi^\e + d^* \psi^\e+\harm\right) \wedge \left(d\phi^\e + d^* \psi^\e+\harm\right)\in\begin{cases}\left[W^{s,q} + L^{\frac{np}{4n-2p}}\right]\left(U; {\rm Ad}\,\eta \otimes \bigwedge^2T^*\R^n\right)\,\, \text{ if } n<p \leq 2n,\\
\left[W^{s,q} + L^{\infty}\right]\left(U; {\rm Ad}\,\eta \otimes \bigwedge^2T^*\R^n\right)\,\, \text{ if } p > 2n
\end{cases}
\end{align*}
\end{small}
by Sobolev embedding.

Repeating the above procedure for finitely many steps, one arrives at
\begin{align}\label{new, x}
\curv - \left(d\phi^\e + d^* \psi^\e+\harm\right) \wedge \left(d\phi^\e + d^* \psi^\e+\harm\right)\in W^{s,q}\left(U; {\rm Ad}\,\eta \otimes \bigwedge^2T^*\R^n\right).
\end{align}
That is, this term is as regular as the prescribed curvature 2-form $\curv$. 

Using once again the elliptic regularity theory, one now concludes from Eq.~\eqref{system for psi epsilon} that $\psi^\e \in W^{s+2,q}\left(U;{\rm Ad}\,\eta \otimes \bigwedge^2T^*\R^n\right)$. Its uniform $W^{s+2,q}$-bound follows from explicit estimates in the Calder\'{o}n--Zygmund theory, and the regularity of $\psi^\e$ up to the boundary $\p U$ follows from the boundary regularity theory of elliptic systems (\emph{cf}. \emph{e.g.}, \cite{adn2}). Hence \emph{Claim~D} is proved.    \end{proof}

\smallskip
\noindent
{\bf Step 6.} As a consequence of Step~5 above, we deduce:

\noindent
{\bf Claim E.} Assume that the prescribed curvature $2$-form $\curv \in W^{s,q}$ for $s \geq 0$, $q \geq \frac{p}{2}>\frac{n}{2}$, and $W^{s,q}\subset L^2$ over $U$, and let $\Omega$ be the connection 1-form in Theorem~\ref{thm: re}. Let $\phi^\e$ be the smooth scalar fields in Step~2, and $\psi^\e$ be the 2-forms with  $W^{s+2,q}$-regularity in Steps~3--5. Define
\begin{equation}\label{new: Omega epsilon def}
\Omega^\e:=d\phi^\e+d^*\psi^\e+\harm \qquad \text{in } U.
\end{equation}
Then, under the smallness condition~\eqref{smallness condition--key}, we have the strong convergence
\begin{align*}
\Omega^\e \longrightarrow \Omega \qquad \text{ in } L^p\left(U;{\rm Ad}\,\eta\otimes\bigwedge^1T^*\R^n\right)
\end{align*}
after passing to subsequences if necessary.

Eq.~\eqref{new: Omega epsilon def} should be compared with Step~1, \emph{Claim~A}, Eq.~\eqref{new, decomposition of Omega}, which motivates the definition of the approximation family $\left\{\Omega^\e\right\}$. 

\begin{proof}[Proof of Claim~E]
By \emph{Claim~D}, $\{d^*\psi^\e\}$ is uniformly bounded in $W^{s,q}\left(U;{\rm Ad}\,\eta\otimes\bigwedge^1T^*\R^n\right)$, which compactly embeds into $W^{1,p}\left(U;{\rm Ad}\,\eta\otimes\bigwedge^1T^*\R^n\right)$. Our goal is to prove that the strong $W^{1,p}$-limit of $\{\psi^\e\}$ equals to the 1-form $\psi$ defined via Eq.~\eqref{psi system} in Step~1, \emph{Claim~A}.

To this end, recall from Step~2, Eq.~\eqref{phi epsilon} that $\phi^\e \to \phi$ in $W^{1,p}$. Then consider the 1-form $$D:=\lim_{\e \to 0} \left(\psi^\e - \psi\right).$$ This limit is \emph{a priori} well defined in $L^p$, thanks to Cauchy--Schwarz and \emph{Claim~C}. (The latter shows that for each $0<\e\leq 1$, the 1-form $\psi^\e$ is uniquely determined under the smallness condition~\eqref{smallness condition--key}; hence there is no ambiguity in defining $D$.) Furthermore, by virtue of  elliptic regularity applied to Eqs.~\eqref{psi system} and \eqref{system for psi epsilon},  $D$ lies in $W^{2, \max\left\{2,\frac{p}{2}\right\}}$ and hence $W^{1,p}$ by Sobolev embedding. (Here, $\psi^\e, \psi\in W^{2,2}$ uses the additional assumption that $W^{s,q} \emb L^2$, so that the prescribed curvature 2-form $\curv \in L^2$.) Repeating the arguments in Step~5, \emph{Claim~D}, we find that $D$ is well defined in $W^{s+2,q}$ over $U$.

Subtracting Eq.~\eqref{psi system} from Eq.~\eqref{system for psi epsilon} and taking $\e \to 0$, we arrive at the following elliptic system for $D$:
\begin{equation}\label{system for D}
\begin{cases}
\Delta D = d^*D \wedge \Omega + \Omega \wedge d^*D + d^*D \wedge d^*D
 \qquad \text{ in } U,\\
\nnn D= 0 \quad \text{ and } \quad \nnn dD=0 \quad \text{ on } \p U,
\end{cases}
\end{equation}
where $\|\Omega\|_{L^p(U)} \leq \kappa$ as in Step~4, \emph{Claim~C}, Eq.~\eqref{smallness condition--key}. Adapting essentially the arguments in Step~4 via exploiting the elliptic estimates for $\T':D \mapsto \Delta^{-1}\left\{d^*D \wedge \Omega + \Omega \wedge d^*D + d^*D \wedge d^*D\right\}$, we find that this map is a coercive mapping from $W^{1,p}$ to $W^{1,p}$. That is, $\left\|\T'D\right\|_{W^{1,p}} \leq (1-\delta)\kappa'$ whenever $\|D\|_{W^{1,p}(U)} \leq \kappa'$, for some $\kappa' = \kappa'(n,p,\kappa)$ and $0<\delta<1$. Thus, Banach fixed point theorem implies that Eq.~\eqref{system for D} has a unique solution if $\|D\|_{W^{1,p}(U)} \leq \kappa'$, which must be $D \equiv 0$. 

Therefore, we have shown that $d^*\psi^\e \to d^*\psi$ in $L^p$, which together with  Step~2, Eq.~\eqref{phi epsilon} implies that $\Omega^\e \to \Omega$ in $L^p$. Hence \emph{Claim~E} is proved.  \end{proof}


\smallskip
\noindent
{\bf Step 7.} 
Define $\mu_U$ to be the smallest positive number such that the following holds:
\begin{align}\label{poincare}
&\left\|\Xi - \frac{1}{{\rm Vol}(U)}\int_U \Xi\,\dd x\right\|_{L^2(U)} \leq \mu_U \Big\{\|d \Xi\|_{L^2(U)}+\|d^* \Xi\|_{L^2(U)}\Big\}\nonumber\\
&\qquad\qquad\text{ for all }  \Xi \in W^{1,2}\left(U;{\rm Ad}\,\eta\otimes\bigwedge^3T^*\R^n\right)\text{ with vanishing normal trace}.
\end{align}
This is equivalent to the usual Poincar\'{e} inequality $$\left\|\Xi-\frac{1}{{\rm Vol}(U)}\int_U \Xi\,\dd x\right\|_{L^2(U)} \leq \mu'_U \|\na \Xi\|_{L^2(U)},$$ where $\na$ is the Riemannian (here, Euclidean) gradient, in view of Gaffney's inequality. 

\noindent
{\bf Claim F.} In the setting as above, and under the smallness condition~\eqref{smallness condition--key} in \emph{Claim~C} (where the constant $\kappa$ depends additionally on $\mu_U$), $d\psi^\e$ is co-closed (\emph{i.e.}, $d^*d\psi^\e=0$) in $U$ for \emph{every} $0<\e\leq 1$.

We remark that the constant $\kappa$ in the $L^p$-smallness condition~\eqref{smallness condition--key} for $\Omega$ is allowed to depend on $U$ through the Poincar\'{e} constant in Eq.~\eqref{poincare}. For Euclidean balls $U = \mathbf{B}_r(a) \Subset \R^n$, one has $\mu_U = Cr$ and, by Eq.~\eqref{eq, kappa} in the proof below, $\kappa \leq C'r^{-1}$, where $C$ and $C'$ are uniform constants independent of $U$.

\begin{proof}[Proof of Claim F]
Recall Eq.~\eqref{system for psi epsilon}, \emph{i.e.}, the elliptic boundary value problems for $\psi^\e$, as well as the definition for $\Omega^\e$ in Step~6, \emph{Claim~E}, Eq.~\eqref{new: Omega epsilon def} and the Hodge decomposition for $\Omega$ in Step~1, \emph{Claim~A}, Eq.~\eqref{new, decomposition of Omega}. As $d\circ\Delta=\Delta\circ d$, we obtain the equation for $d\psi^\e$ by taking $d$ to Eq.~\eqref{system for psi epsilon} and utilising Eqs.~\eqref{curvature, def}, \eqref{psi system}:
\begin{align}\label{new, dpsi eq}
\Delta d\psi^\e &= d(\Omega \wedge \Omega) - d\left(\Omega^\e \wedge \Omega^\e\right) \nonumber\\
&= d\Omega \wedge \Omega -\Omega \wedge d\Omega - d\Omega^\e \wedge \Omega^\e + \Omega^\e \wedge d\Omega^\e\nonumber\\
&= dd^*\psi \wedge \Omega - \Omega \wedge dd^*\psi - dd^*\psi^\e \wedge \Omega^\e + \Omega^\e \wedge dd^*\psi^\e.
\end{align}
The final line holds by $d\Omega =d(d\phi+d^*\psi+\harm)=dd^*\psi$ and, similarly, $d\Omega^\e = dd^*\psi^\e$. This equation is supplemented with the boundary condition $\nnn \left(d\psi^\e\right)=0$ on $\p U$.

We note in passing that we have an ``incomplete'' boundary condition here, in the sense that $\nnn \left(d\psi^\e\right)=0$ does not make Eq.~\eqref{new, dpsi eq} into a elliptic boundary value problem of PDE systems in the sense of Agmon--Douglis--Nirenberg \cite{adn1, adn2}. Nonetheless, in the arguments below ellipticity plays no role; all we need is integration by parts using this incomplete boundary condition.

By virtue of Step~5, \emph{Claim~D}, one has that $\psi^\e \in W^{2,2}\left(U; {\rm Ad}\,\eta \otimes \bigwedge^2T^*\R^n\right)$. Hence,  Eq.~\eqref{new, dpsi eq} is well defined in the sense of distributions, and so is the $L^2$-inner product $\bra d\psi^\e, \Delta d\psi^\e\ket_U$. An application of the Stokes' theorem (or the Gauss--Green theorem) yields that
\begin{align}\label{new, 1}
\bra d\psi^\e, \Delta d\psi^\e\ket_U&=\|d^* d\psi^\e\|_{L^2(U)}^2 + \int_{\p U} \ttt\left(d^*d\psi^\e\right)\wedge\star \nnn(d\psi^\e)\nonumber\\
&= \|d^*d\psi^\e\|_{L^2(U)}^2,
\end{align}
thanks to the boundary condition $\nnn \left(d\psi^\e\right)=0$ on $\p U$.

On the other hand, Eq.~\eqref{new, dpsi eq} and H\"{o}lder's inequality give us
\begin{align*}
\bra d\psi^\e, \Delta d\psi^\e\ket_U \leq C_3\Big\{\|\Omega^\e\|_{L^p(U)}\|dd^*\psi^\e\|_{L^2(U)} + \|\Omega\|_{L^p(U)}\|dd^*\psi\|_{L^2(U)}  \Big\} \|d\psi^\e\|_{L^\frac{2p}{p-2}(U)} 
\end{align*}
where $C_3=C_3(p,n)$. Hence, by Eq.~\eqref{new, 1} and Step~6, \emph{Claim~E},
\begin{align*}
\|d^*d\psi^\e\|_{L^2(U)} \leq 2C_3\Big\{\|\Omega\|_{L^p(U)}\|dd^*\psi^\e\|_{L^2(U)} + \mathfrak{o}(\e) \Big\} \|d\psi^\e\|_{L^\frac{2p}{p-2}(U)} 
\end{align*}
for some $\mathfrak{o}(\e)\to 0$ as $\e \to 0^+$. Recall from Step~3 that $$
dd^*\psi^\e = d\Omega^\e = \curv - \Omega^\e\wedge \Omega^\e$$ and from Step~5, Eq.~\eqref{new, x} that $$\left\|\curv - \Omega^\e\wedge \Omega^\e\right\|_{L^2(U)} \leq C_4(n).$$ In view of the $L^p$-smallness of $\Omega$ (\emph{Claim~C}, Eq.~\eqref{smallness condition--key}) and the compact embedding $
W^{1,2}\emb\emb L^{\frac{2p}{p-2}}$ over $U \Subset \R^n$ for $p>n$, one thus has
\begin{align*}
\|d^*d\psi^\e\|_{L^2(U)} \leq 2C_3\Big\{ C_4\kappa + \mathfrak{o}(\e) \Big\}\Big\{\delta \|d^*d\psi^\e\|_{L^2(U)} + C_\delta \left\|d\psi^\e\right\|_{L^2} \Big\},
\end{align*}
where $\delta>0$ is arbitrarily small and $C_\delta \sim \delta^{-1}$.

Without loss of generality we assume $\mathfrak{o}(\e) \lesssim \kappa \leq 1$, so that
\begin{align*}
\|d^*d\psi^\e\|_{L^2(U)} \leq C_5\kappa\Big\{\delta \|d^*d\psi^\e\|_{L^2(U)} + C_\delta \left\|d\psi^\e\right\|_{L^2} \Big\}
\end{align*}
for some $C_5=C_5(p,n)$. Now, take first $$\delta = \frac{1}{2C_5}$$ and then 
\begin{equation}\label{eq, kappa}
\kappa \leq \min\left\{1,\,\frac{1}{4C_5C_\delta \mu_U}\right\}
\end{equation}
 to ensure that 
\begin{align*}
\|d^*d\psi^\e\|_{L^2(U)} \leq \frac{1}{2\mu_U} \|d\psi^\e\|_{L^2(U)},
\end{align*}
where $\mu_U$ is the Poincar\'{e} constant in Eq.~\eqref{poincare} at the beginning of this step. 
Taking $\Xi = d\psi^\e$ in Eq.~\eqref{poincare} and utilising the identities $dd\psi^\e=0$ and $\int_U d\psi^\e\,\dd x = 0$ for $\nnn\psi^\e=0$ on $\p U$ (the latter one follows from the Stokes' theorem), we infer that
\begin{align*}
\|d^*d\psi^\e\|_{L^2(U)} \leq \frac{1}{2\mu_U} \|d\psi^\e\|_{L^2(U)} \leq \frac{1}{2}\|d^*d\psi^\e\|_{L^2(U)}.
\end{align*} 
Hence $d^*d\psi^\e \equiv 0$, which proves \emph{Claim~F}.  \end{proof}

\smallskip
\noindent
{\bf Step 8.} Now we check that $\left\{\Omega^\e\right\}_{0<\e\leq 1}$ is indeed a smooth approximating family of $\Omega$ with the same prescribed curvature 2-form; \emph{i.e.}, for every $0<\e \leq 1$, it holds that
\begin{equation}\label{new, z, eq}
d\Omega^\e + \Omega^\e \wedge \Omega^\e = \curv.
\end{equation}

For this purpose, note by Step~3, Eq.~\eqref{system for psi epsilon} that $\Delta \psi^\e = \curv- \Omega^\e \wedge \Omega^\e$, where $\Omega^\e = d\phi^\e + d^*\psi^\e + \harm$ as in  Step~6, Eq.~\eqref{new: Omega epsilon def}. So we just need to show that 
\begin{equation}\label{new, y}
\Delta\psi^\e = d\Omega^\e.
\end{equation}
But Eq.~\eqref{new: Omega epsilon def} gives us $$d\Omega^\e = dd^*\psi^\e = \Delta \psi^\e - d^*d\psi^\e,$$ whereas Step~7, \emph{Claim~F} yields that $
d^*d\psi^\e \equiv 0$ in $U$. Now Eq.~\eqref{new, y} and hence Eq.~\eqref{new, z, eq} follows.

\smallskip
\noindent
{\bf Step 9.} We first conclude the proof for Theorem~\ref{thm: local}. One only needs to check that the arguments in the previous Steps~1--8 carry over to $U \Subset \M$ for a general Riemannian manifold $\M$. This is because all the theories of elliptic PDE systems, Hodge decomposition, and Sobolev embeddings used in Steps~1--8  extend to general $\M$ (\cite{sch, h}). In fact, the arguments in Steps~1--8 have been formulated so that they carry over to bounded smooth domains of Riemannian manifolds almost verbatim, once all occurrences of $\R^n$ are replaced by $\M$, and all the constants $C_i$ and $\kappa_i$ in Steps~1--8, which depend only on $p$ and $n$ therein, are now chosen to depend only on $p$ and the geometry of $\M$.

The only slightly more involved constant is $\kappa$, the upper bound in the $L^p$-smallness condition~\eqref{smallness condition--key} for $\Omega$. For this purpose, in view of Step~7, Eq.~\eqref{eq, kappa}, the constant $\kappa$ now depends on $C_5$ (which in turn depends only on $p$ and $\M$ by the last paragraph) and the Poincar\'{e} constant $\mu_U$. Note, however, that if $\mu_U$ can be chosen uniformly bounded from the above, then by Eq.~\eqref{eq, kappa} we can select $\kappa=\kappa(p,\M)$ independently of $U$. This is the case for compact $\M$, where $\mu_U$ is no larger than a uniform constant $C=C(\M)$ times the  radius of $U$ in the Riemannian metric. See, \emph{e.g.}, the recent work \cite{gis} by  Grigor'yan--Ishiwata--Saloff-Coste  for a survey of this result and  far-reaching generalisations to open manifolds.

Thus, the proof for Theorem~\ref{thm: local} is now complete.

Finally, consider a closed (\emph{i.e.}, compact, boundaryless) Riemannian manifold $\M$. To conclude Theorem~\ref{thm: re}, we replace in Steps~1--8 above all the occurrences of subdomains $U \Subset \M$ by the whole manifold $\M$. In this regard, all the boundary conditions imposed on $\p U$ should be dropped, since $\p \M = \emptyset$.  For instance, Step~1, \emph{Claim~A}, Eq.~\eqref{psi system} should be substituted by $\Delta \psi = \curv - \Omega \wedge \Omega$ over $\M$. The existence and regularity for $\psi$ remain valid (see \cite{sch, h}), while the uniqueness holds modulo the addition of any harmonic 1-form. Similarly for Step~3, Eq.~\eqref{system for psi epsilon}, \emph{i.e.}, the PDE for $\psi^\e$.

However, the absence of boundary conditions imposes no severe difficulty in our case, as in each of the relevant places we may fix once and for all a harmonic form on $\M$ to retain the uniqueness of solution. In this way, all the arguments in  Steps~1--8 carry out verbatim, once we replace each occurrence of $U$ and $\R^n$ with $\M$. In particular, our smoothing construction for $\Omega$ is completely global: the smoothing sequence $\{\Omega^\e\}$ is constructed via global solutions to several elliptic PDE systems over $\M$, rather than obtained from glueing/patching together smoothing constructions in local charts.

This completes the proof for Theorem~\ref{thm: re}.   \end{proof}

\begin{remark}
It would be interesting to see if Theorem~\ref{thm: main} (or equivalently, Theorem~\ref{thm: re}) hold for the endpoint case $p=n$ under the assumption that $\|\Omega\|_{L^n}$ is small. For this case, probably clever applications of gauge transforms are necessary, as in Rivi\`{e}re's seminal work \cite{riviere}.
\end{remark}

\section{Direct consequences of Theorem~\ref{thm: main}}\label{sec: existence}

Here we present the proofs for two direct consequences of the smoothability Theorem~\ref{thm: main}. These are Corollary~\ref{cor: regularity}, the ``regularity pending uniqueness'' result, and Proposition~\ref{propn: luzin}, the Luzin-type smoothability theorem without $L^p$-smallness condition~\eqref{smallness condition--key} on $\Omega$. The arguments for Proposition~\ref{propn: luzin} contain a recurring theme, which will reappear in the proof of Theorem~\ref{thm: luzin, almost global isom imm}.

\begin{proof}[Proof of Corollary~\ref{cor: regularity}]

Theorem~\ref{thm: main} yields the existence of smooth solutions $\left\{\Omega^\e\right\} \subset \conn(\M,\kappa_0,\curv)$ to the equation $d\Omega^\e + \Omega^\e \wedge \Omega^\e=\curv$ such that $\Omega^\e \to \Omega \in  \conn^p_0(\M,\kappa_0,\curv)$ strongly in $L^p$. Under the assumption of uniqueness,  
$\Omega^\e \equiv \Omega$ for every $\e$.   \end{proof}

\begin{remark}\label{rem: nonunique}
In general, one should not expect uniqueness of the connection 1-form, even if the associated curvature 2-form is $C^\infty$ and small (\emph{e.g.}, in the $W^{-1,p}$-topology). Indeed, in the case of flat curvature, namely that $\curv\equiv 0$, Wehrheim \cite{w-new} showed that any $L^p$ connection $\Omega$ with $p>n$ is gauge-equivalent to a $C^\infty$-connection. That is, there is a $W^{2,p}$-gauge transform $u$ such that $u^\#\Omega$ is $C^\infty$. But then $u^\#\Omega$ is also flat, since the curvature 2-form transforms under the action of gauges via $F_{u^\#\Omega} = u^{-1} F_\Omega u$.
\end{remark}

\begin{proof}[Proof of Proposition~\ref{propn: luzin}]

Let $\mathcal{A} = \{\U_i\}_{i \in \mathcal{I}}$ be an atlas for $(\M,g)$ and $\{\chi_i\}_{i \in \mathcal{I}}$ be a $C^\infty$-partition of unity subordinate to this atlas. We choose $\mathcal{A}$ such that the volume of each $\mathcal{U}_i$ is less than $v_0$. Hence, by assumption, there are local connection 1-forms $\Omega_i \in \conn^p_0(\U_i)$ such that $$d \Omega_i + \Omega_i\wedge \Omega_i=\curv\big|_{\U_i}.$$ 
Then set
 \begin{align*}
 \tilo := \sum_{i \in \mathcal{I}} \chi_i \Omega_i.
 \end{align*}
This is a well-defined element in $\conn^p_0(\M)$: connection 1-forms take value in the vector space $\lie \otimes \G(T^*\M)$, so the naive glueing construction via partition of unity works here\footnote{This is not the case for the construction of global \emph{gauges} in Uhlenbeck's works \cite{u', u}: gauge transforms, by definition, take values in Lie groups, so the naive linear combination ceases to be valid.}.

Note that
\begin{align*}
d\tilo + \tilo \wedge \tilo &=  \sum_{i \in \mathcal{I}} \left(d\chi_i \wedge \Omega_i + \chi_i \,d\Omega_i\right) + \sum_{(i,j) \in \mathcal{I}\times \mathcal{I}}\chi_i\chi_j \Omega_i \wedge \Omega_j\\
&=  \sum_{i \in \mathcal{I}}d\chi_i \wedge \Omega_i + \sum_{i \neq j}\chi_i\chi_j \Omega_i \wedge \Omega_j + \sum_{i \in \mathcal{I}} \chi_i \left\{d\Omega_i + \chi_i\Omega_i \wedge \Omega_i\right\}\\
&=  \sum_{i \in \mathcal{I}}d\chi_i \wedge \Omega_i + \sum_{i \neq j}\chi_i\chi_j \Omega_i \wedge \Omega_j + \curv -\sum_{i \in \mathcal{I}} \chi_i\left(1-\chi_i\right)\Omega_i \wedge \Omega_i.
\end{align*}
The final line follows from $d \Omega_i + \Omega_i\wedge \Omega_i=\curv\big|_{\U_i}$ and the construction of $\{\chi_i\}$.

Now let $$\M_\circ:= \left\{x \in \M:\,d\tilo + \tilo \wedge \tilo=\curv\right\}$$ be the ``good'' set. Then
\begin{align*}
\M\sim\M_\circ \subset \bigcup_{i\in\mathcal{I}} \left\{d\chi_i = 0\right\}  \cup   \bigcup_{i\in\mathcal{I}} \left\{\chi_i(1-\chi_i) = 0\right\} \cup \bigcup_{\left\{(i,j) \in \mathcal{I}\times \mathcal{I}:\, i \neq j\right\}} \left\{\chi_i\chi_j = 0\right\}.
\end{align*}
The set of the right-hand side can be made to have Riemannian volume less than $\e$: we first choose the atlas $\mathcal{A}$ --- upon further refinements if necessary --- such that each chart (which is already known to have volume $<v_0$) intersects all the other charts in a set of arbitrarily small total volume, and then choose the subordinating partition of unity such that each $\chi_i$ transits from $1$ to $0$ in a subset of $\U_i$ of arbitrarily small volume. 

The proof is now complete, once we notice that $\tilo$ has the same regularity as each $\Omega_i$. 
\end{proof}

The proof of Proposition~\ref{propn: luzin} carries over modulo obvious modifications for noncompact $\M$ admitting a ``tame'' $C^\infty$-partition of unity, in the sense that there is a sufficiently refined atlas such that each chart intersects $\leq N$ other charts, where $N$ is a uniform  number. This is the case for a wide class of Riemannian manifolds, \emph{e.g.}, for an open (namely, complete and noncompact) manifold $\M$ with bounded geometry. Remark~\ref{remark, ppl G bdl}  remains valid in this setting.

\section{Applications to Isometric Immersions}\label{sec: isom imm}

Now let us explore the applications of Theorem~\ref{thm: main} to the isometric immersion problem. Once again, we refer  to Han--Hong \cite{hh} for a comprehensive survey on this classical problem.


\subsection{Fundamental theorem of submanifold theory with lower regularity}

This subsection is mainly dedicated to the proof of Theorem~\ref{thm: existence}, which is a $W^{2,p}$-version of the ``fundamental theorem of submanifold theory''. This terminology is coined as an extension of the ``fundamental theorem of surface theory'' which, under suitable regularity assumptions, ascertains the equivalence between the local existence and uniqueness modulo Euclidean rigid motions of isometric immersions and the validity of Gauss--Codazzi--Ricci equations.

As  remarked in \S\ref{sec: Intro}, our proof is essentially a recast of the arguments in \cite[\S 5]{cl}, which utilises ideas from Tenenblat \cite{ten} and S. Mardare \cite{m05, m07} as well as Theorem~\ref{thm: main}. Note, however, that Theorem~\ref{thm: existence} (2) is a new result, which establishes the existence of global $W^{2,p}$-isometric immersions for ``nearly flat'' $(\M,g)$ given the Gauss--Codazzi--Ricci equations without any topological restrictions on $\M$. Meanwhile, we take this opportunity to elucidate some issues on local \emph{vs.} global isometric immersions in the literature.

We begin with a quick review of Theorem~\ref{thm: existence, early}.

\begin{proof}[Sketched proof of Theorem~\ref{thm: existence, early}]

It is well known that (II) is equivalent to (III). Also, it is a standard result that (III) is a necessary condition for (I) (\cite[Chapter~6]{d}). All the above hold as purely algebraic  identities, which can be easily validated in the sense of distributions.

Now we turn to (II) $\Rightarrow$ (I). We assume the Cartan formalism --- in particular, the second structural equation $d\Omega + \Omega \wedge \Omega=0$ in the sense of  distributions, where 
\begin{equation}\label{Omega, new'}
\Omega = \begin{bmatrix}
\na & \two\\
-\two & \na^\perp
\end{bmatrix} \in L^p_\loc(\M;T^*\M \otimes \mathfrak{so}(n+k,\R))
\end{equation}
is the connection 1-form defined via Eqs.~\eqref{Omega 1}--\eqref{Omega 3}. This can also be regarded as the equation for flat curvature; see Eq.~\eqref{curvature, def}. What we seek is a \emph{local} isometric immersion $\iota: \U\subset \left(\M^n,g\right) \emb \R^{n+k}$.  As in the proof for  \cite[Theorem~5.2]{cl}, which is modelled over the classical arguments in \cite{ten}, the isometric immersion  can be obtained by solving first a \emph{Pfaff system}
\begin{align}\label{pfaff, isom imm}
d\wp= -\Omega \wp
\end{align}
and then a \emph{Poincar\'{e} system}
\begin{equation}\label{poincare, isom imm}
d\iota = \omega \wp.
\end{equation}
Here, $\omega \in \G(T^*\M \oplus \{{\bf 0}\})$ (where ${\bf 0}$ is the zero-section of $E$) is given by 
\begin{equation*}
\omega:=\big(\omega^1,\ldots,\omega^n,0,\ldots,0\big),
\end{equation*}
and $\wp \in W^{1,p}(\M;T^*\M \oplus E)$ is an intermediate variable to be solved. 

Now we take local charts $\U \Subset \U'\Subset \U''\Subset \M$, where $\U''$ is sufficiently small such that $E$ is trivialised over it, and that $\left\|\Omega\big|_{\U''}\right\|_{L^p}<\kappa_0$. Here $\kappa_0$ is the small constant as in Theorem~\ref{thm: main}; that is, we assume that $\Omega\big|_{\U''}$ lies in the moduli space $\conn^p_0\left(\M,\kappa_0,\curv \equiv 0\right)$. Here, our choice of the bundle and structure group are $\eta = \T^*\M \oplus E$ and $G=SO(n+k;\R)$, respectively. Also fix a base point $x_\star \in \U$ and data $\wp_\star \in T^*_{x_\star}\M \oplus E\big|_{x_\star}$ and $\iota_\star \in \R^{n+k}$.

By assumption $d\Omega\big|_{\U''}+\Omega\big|_{\U''}\wedge\Omega\big|_{\U''}=\curv$, so we can apply Theorem~\ref{thm: main} to find a smooth family of connection 1-forms over $\U''$ with the prescribed curvature 2-form: $$\left\{\Omega^\e\right\}_{0<\e\leq 1} \subset \conn(\U; \T^*\M \oplus E) \cap \conn^p_0\left(\M,\kappa_0,\curv \equiv 0\right),$$ such that $\Omega^\e \to \Omega$ in $L^p$ on $\U''$, thanks to the $L^p$-smallness of $\Omega\big|_{\U''}$. Classical arguments in \cite{ten} via the Frobenius theorem (see also \cite[Step~5]{cl}) established that\footnote{Let us point out that slightly different arguments can be found in S. Mardare \cite{m05, m07} and Szopos \cite{sz}, which also applied the Frobenius theorem on integrability to locally solve for Pfaff and Poincar\'{e} systems, while utilising various matrices or tensors that probably motivated by the literature in elasticity theory. In contrast, Tenenblat in \cite{ten} worked with more ``geometrically natural'' quantities in the framework of the method of moving frames. Some of the ideas in \cite{m05, m07, sz} pertaining to the fundamental theorem of submanifold theory (with low regularity) may be connected to those in Ciarlet--Larsonneur \cite{ciarlet'}. Also, it is commented in \cite{ciarlet'} that demonstrations of the \emph{local} version of the fundamental theorem of surface theory (\emph{i.e.}, $n=2$ and $k=1$) can be found in Choquet-Bruhat, Dewitt
Morette, and Dillard-Bleick \cite[p. 303]{add} and Malliavin \cite[p. 133]{add'}, while the outline for the proof of the global version (with $C^3$-isometric immersions) in \cite{ciarlet'} follows that of the elasticity paper by Blume \cite{b}. We also note that arguments and computations in  S. Mardare \cite{m05} also play a key role in the extension of the fundamental theorem of surface theory to $W^{2,2}$-isometric immersions. See Litzinger \cite{l}.}
\begin{itemize}
\item
The Pfaff system~\eqref{pfaff, isom imm} $d\wp^\e = -\Omega^\e \wp^\e$ subject to $\wp^\e(x_\star)=\wp_\star$ has a unique smooth solution $\wp^\e$ over $\U'$ for each $0<\e\leq 1$, such that $\Omega^\e \mapsto \wp^\e$ is $L^p\to W^{1,p}$ continuous.
\item
The Poincar\'{e} system~\eqref{poincare, isom imm} $d\iota^\e = \omega \wp^\e$ subject to $\iota^\e(x_\star)=\iota_\star$ has a unique smooth solution $\iota^\e$ over $\U$ for each $0<\e\leq 1$, such that $\wp^\e\mapsto \iota^\e$ is $W^{1,p}\to W^{2,p}$ continuous.
\end{itemize}
In this way, by sending $\e \to 0^+$ one obtains $\iota^\e \to \iota$ in $W^{2,p}$, and it is direct to check that $\iota$ is an isometric immersion on $\U$ with associated connection 1-form coinciding with $\Omega$. 

Finally, the uniqueness of $\iota$ modulo Euclidean rigid motions in $\R^{n+k}$ follows from the uniqueness of solutions to the Pfaff and Poincar\'{e} systems, as well as the freedom of choosing orthonormal frames on $\R^{n+k}$ (in other words, translations and/or rotations of isometric immersions remain isometric immersions).  \end{proof}

In Theorem~\ref{thm: existence, early} we have assumed that $\curv = 0$, \emph{i.e.}, the connection 1-form in Eq.~\eqref{Omega, new'} associated to the Cartan formalism is flat. The case $\curv \neq 0$ corresponds to isometric immersions of $(\M,g)$ into a general ambient manifold $(\mathcal{N},h)$ that is not necessarily Euclidean. In this case, the Gauss--Codazzi--Ricci equations (and hence the Cartan formalism) are still well-known.  See \cite[Chapter~6]{d}. But the existence of isometric immersions assuming the Gauss--Codazzi--Ricci equations --- even the local version with $C^\infty$-regularity --- seems to be missing from the literature. We plan to address this issue in follow-up works.

Now we turn to the proof of Theorem~\ref{thm: existence}, which is a global version of Theorem~\ref{thm: existence, early} above.

\begin{proof}[Proof of Theorem~\ref{thm: existence}]


We first observe that if $\M=\M^n$ is a ``nearly flat'' closed Riemannian manifold as in  Theorem~\ref{thm: existence} (2), then the $L^p$-connection 1-form $\Omega$ (see Eq.~\eqref{Omega, new'}, with $\curv_\Omega\equiv 0$ due to the isometric immersion into Euclidean space, and $p>n$) can be \emph{globally} smoothed by Theorem~\ref{thm: main}. That is, we can find a family of $C^\infty$-connection 1-forms $\left\{\Omega^\e\right\} \subset \conn\left(\M,\kappa_0,\curv\equiv 0\right)$ such that $\Omega^\e \to \Omega$ strongly in $L^p$. From here, the arguments based on Pfaff and Poincar\'{e} systems in the proof of Theorem~\ref{thm: existence, early} carry over.

Now we turn to the case that $\M$ is simply-connected. Let us first prove for $\M$ compact. Take an atlas of open charts $\{\U_i\}_{i \in \mathcal{I}}$ for $\M$ such that
\begin{itemize}
\item
the indexing set $I$ is finite;
\item
each $\U_i$ is contractible;
\item
 on each $\U_i$ the bundle $E$ is trivialised; and
 \item
 each $\U_i$ is so small that $
\|\Omega\|_{L^p(\U_i)} \leq \kappa_0/2$.
\end{itemize} Here $\kappa_0 = \kappa_0(p,\M)$ is the small constant as in Theorem~\ref{thm: main} such that $\conn^p_0(\U_i,\kappa_0,\curv\equiv 0)$ is smoothable for all $i \in \mathcal{I}$. In this case, although we do not have a global smoothing of $\Omega$ over the whole manifold $\M$, the local isometric immersions arising from the local smoothing can be easily glued together to form a global $\iota: (\M,g) \emb \R^{n+k}$.


To be more precise, let us take  atlases $\{\V_i\}$, $\{\V_i'\}$,  and $\{\V_i''\}$ for $\M$ such that each $\V_i \Subset \V_i' \Subset  \V_i'' \Subset \U_i$ and $\bigcup_{i \in \mathcal{I}}\V_i = \M$. As in the proof of Theorem~\ref{thm: existence, early}, we may approximate $\Omega\big|_{\V_i''}$ by smooth connection 1-forms $\left\{\Omega^\e_i\right\}\subset \conn\left(\V_i'', \kappa_0, 0\right)$ on each $\V_i''$, solve for $\wp_i^\e$ from the Pfaff system \begin{equation}\label{loc pfaff}
d\wp_i^\e = -\Omega_i^\e\wp_i^\e\qquad\text{ on } \V_i',
\end{equation}
and then solve for $\iota_i^\e$ from the Poincar\'{e} system
\begin{align}\label{loc poincare}
d\iota_i^\e = \omega \wp_i^\e\qquad \text{ on } \V_i.
\end{align}

The previous arguments can be globalised via a straightforward glueing construction, again assuming that $\M$ is compact, by induction on $i \in \mathcal{I}$ = the indexing set for the atlases. Without loss of generality, say $\mathcal{I}= \{1,\ldots, N\}$. Relabelling the indices if necessary, we may assume that for each $L \in \{2,\ldots,N\}$ it holds that
\begin{equation}\label{chart condition}
\V_L \cap \left( \bigcup_{j=1}^{L-1}\V_j \right) \neq \emptyset.
\end{equation}

\begin{itemize}
\item
For $i=1$, fix an arbitrary point $x_\star \in \V_1$ and arbitrary values $\xi_1 \in T^*_{x_\star}\M\otimes E\big|_{x_\star}$ and $a_1 \in \R^{n+k}$. We solve the Pfaff system~\eqref{loc pfaff} on $\V_1'$ under the condition $\wp_1^\e(x_\star)=\xi_1$, and solve the Poincar\'{e} system~\eqref{loc poincare} on $\V_1$ under the condition $\iota^\e_1(x_\star) = a_1$. Both $\wp_1^\e$ and $\iota^\e_1$ are unique in $\V_1$ by previous arguments. 

\item
Suppose that we have constructed for some $K \in \{1,2,\ldots,N-1\}$ the maps $\wp^\e \in C^\infty\left(\bigcup_{i=1}^K \V_i'; T^*\M\otimes E\right)$ and $\iota^\e \in C^\infty\left(\bigcup_{i=1}^K \V_i;\R^{n+k}\right)$  such that, for each $i \in \{1,\ldots,K\}$, the restricted maps $\wp^\e_i := \wp^\e\big|_{\V_i'}$ and $\iota^\e_i := \iota^\e\big|_{\V_i}$ satisfy the Pfaff system~\eqref{loc pfaff} on $\V_i
'$ and the Poincar\'{e} system~\eqref{loc poincare} on $\V_i$, respectively. The goal is to extend the definitions of $\wp^\e$ to $\V_{K+1}'$ and $\iota^\e$ to $\V_{K+1}$, respectively.

Let us first consider $\wp^\e$. For this purpose, the condition~\eqref{chart condition} ensures that $\V_{K+1}'$ intersects some chart $\V_j'$ nontrivially, where $j \in \{1, \ldots, K\}$. Pick $x_{j,\star} \in \V_j' \cap \V_{K+1}'$ and solve the  Pfaff system~\eqref{loc pfaff} on $\V_{K+1}'$, subject to the condition $\wp_{K+1}^\e(x_{j,\star}) = \wp^\e(x_{j,\star})$; here the right-hand side has already been given in the earlier steps of the induction, and is equal to $\wp^\e_j(x_{j,\star})$. By the existence and uniqueness of solution to the Pfaff system, the solution $\wp^\e_{K+1}$ exists on $\V_{K+1}'$ and coincides with $\wp^\e_j$ on $\V_{K+1}'\cap\V_j'$. Then, for any other $\tilde{j} \in  \{1, \ldots, K\}$ with $\V_{K+1}'\cap\V_{\tilde{j}}' \neq \emptyset$, by virtue of condition~\eqref{chart condition} one can find a subset $\{j_1, \ldots, j_\ell\}$ of $\{1,\ldots,K\} \sim \left\{j,\tilde{j}\right\}$ such that $\V_j' \cap \V_{j_1}'$, $\V_{j_1}' \cap \V_{j_2}'$, $\ldots$, $\V_{j_{\ell-1}}' \cap \V_{j_{\ell}}'$, and $\V_{j_\ell}' \cap \V_{\tilde{j}}'$ are all nonempty. But $\wp^\e$ is already known to be well-defined on $\V_1' \cup \ldots \cup \V_K'$ by  induction hypothesis. Then $\wp^\e_{K+1}$ coincides with $\wp^\e_{\tilde{j}}$ on $\V_{K+1}'\cap\V_{\tilde{j}}'$ too, in view of the assumption that $\M$ is simply-connected$^{{(\dagger)}}$. In this way, we succeed in extending $\wp^\e$ to the domain $\bigcup_{i=1}^{K+1} \V_i'$. This completes the induction.  
\end{itemize}

In the above inductive argument, the key step that essentially relies on the simple connectedness of $\M$ is ${{(\dagger)}}$. It is equivalent to the following statement: 
\begin{quote}
Fix a point $x \in \M$ and any loop $\gamma \subset \M$ based at $x$. Cover $\gamma$ by a chain of charts $U_1, U_2, \ldots, U_L$ such that $x \in U_1 \cap U_L$ and each $U_j$, $1 \leq j \leq L$ intersects only with $U_{j-1}$ and $U_{j+1}$ (with $j$ taken mod $L$). The two $C^\infty$-solutions to the Pfaff system~\eqref{loc pfaff}  coincide at $x$ --- one is obtained by solving the equations only in $U_1$, and the other is obtained by solving in $U_1, U_2, \ldots, U_L$ in this order, via the extension process as above.
\end{quote}
To see this, by solving Eqs.~\eqref{loc pfaff} and \eqref{loc poincare} we obtain smooth isometric immersions $\iota_i^\e$ on those charts $U_1, U_2, \ldots, U_L$  (in this order) along $\gamma$. But $\M$ is simply-connected, so the values of $d\iota_1^\e$ and $d\iota_L^\e$ coincide on $U_1 \cap U_L$. In view of the Poincar\'{e} system~\eqref{loc poincare}, we thus infer that the solutions to the Pfaff system~\eqref{loc pfaff} coincide, \emph{i.e.},  $\wp_1^\e=\wp_L^\e$ on $U_1 \cap U_L$\footnote{In passing, note that the situation is different for the Poincar\'{e} system~\eqref{loc poincare}: we do not claim $\iota^\e_1 = \iota^\e_L$  on $U_1 \cap U_L$. This subtle difference arises because we are considering isometric immersions rather than embeddings.}.

From the previous inductive arguments, in the case of  compact $\M$, we have obtained a \emph{globally} defined map $$\wp^\e \in C^\infty(\M; T^*\M \otimes E)$$ such that the restrictions $\wp^\e_i \equiv  \wp^\e\big|_{\V_i'}$ satisfy the Pfaff system~\eqref{loc pfaff} for each $i \in  \mathcal{I}$. From this we can also solve the Poincar\'{e} system~\eqref{poincare, isom imm}  \emph{globally} as in \cite{ten}, and the global solution $\iota^\e$ must coincide with the solutions $\iota^\e_i$ to the (already obtained) local Poincar\'{e} system~\eqref{loc poincare} when restricted to each $\V_i$, modulo a Euclidean rigid motion applied to each solution over $\V_i$.

As in the proof of Theorem~\ref{thm: existence, early}, the maps $\Omega^\e_i \mapsto \wp^{\e}_i$ and $\wp^\e_i \mapsto \iota_i^\e$ are $L^p \to W^{1,p}$ continuous over $\V'_i$ and $W^{1,p} \to W^{2,p}$ continuous over $\V_i$, respectively, for each $i \in \mathcal{I}$. Sending $\e \to 0^+$ and using the definition of Sobolev spaces on compact manifolds (\cite{h}), one arrives at a global isometric immersion $\iota \in W^{2,p}(\M,\R^{n+k})$. This completes the proof for compact $\M$.

Finally we turn to the case that $\M$ is noncompact. To this end, take a compact exhaustion\footnote{The existence of compact exhaustion on Riemannian manifold is a standard result in topology. Indeed, it is known that on a topological space which is Hausdorff, locally Euclidean, and connected, the following are equivalent: (i), second countability; (ii), paracompactness; and (iii), the existence of compact exhaustion.} $\M_1 \Subset \M_2 \Subset \M_3 \Subset \ldots \nearrow \M$. The earlier parts of the proof yield a $W^{2,p}$-isometric immersion $\iota$ for $\M_2$. We can then utilise the glueing arguments as above to extend it to $\M_3$, by considering a suitable atlas for ${\rm int}\left(\M_3\right) \sim \M_1$. Again, there is no obstruction for such glueing over simply-connected $\M$.  Continuing in this way, we may extend the isometric immersion $\iota$ (not relabelled) in countably many steps to the whole manifold $\M$, whose restriction to each $\M_n$ (hence to each compact subset of $\M$) is in $W^{2,p}$. This shows that the isometric immersion $\iota$ is in $W^{2,p}_\loc(\M,\R^{n+k})$, which completes the proof of Theorem~\ref{thm: existence} (1).  \end{proof}

The ``almost global'' version of this result in the Luzin sense, namely, Theorem~\ref{thm: luzin, almost global isom imm}, is a straightforward consequence.

\begin{proof}[Proof of Theorem~\ref{thm: luzin, almost global isom imm}]

Fix a triangulation of $\M$, and let $\M_\times$ be a smooth closed neighbourhood of the 1-skeleton of $\M$ with arbitrarily small Riemannian volume. We apply Theorem~\ref{thm: existence} to the simply-connected submanifold $\M \sim \M_\times$ to conclude.  \end{proof}

Finally we have a weak compactness theorem for immersed Euclidean submanifolds of regularity $W^{2,p}$ for $p>n=\dim\M$. Here 
$\delta$ is the Euclidean metric on $\R^{n+k}$, and $\na^{\perp,\ell}$ is the orthogonal projection of the Levi-Civita connection on $\R^{n+k}$ to the normal bundle of $f^\ell$.  Our formulation follows  \cite[Theorem~4]{l}; see also Ciarlet--C. Mardare \cite{new1} and Chen--Li--Slemrod \cite{new2}. 

\begin{corollary}\label{cor: weak compactness, p>n}
Let $\M$ be an $n$-dimensional connected and simply-connected  closed manifold. Let $\left\{f^\ell\right\} \subset W^{2,p}(\M, \R^{n+k})$ be a sequence of immersions with uniformly bounded $W^{2,p}$-norms, $p>n$, with corresponding metrics $\left\{g^\ell := \left(f^\ell\right)^\# \delta\right\}$, second fundamental forms $\left\{\two^\ell\right\}$, and affine normal connections $\left\{\na^{\perp,\ell}\right\}$. Assume in addition that $\left\{g^\ell\right\}$ is \emph{non-degenerate}, in the sense that   $$\inf_{x \in \M,\,\ell\in\mathbb{N}} \Big\{\text{the lowest eigenvalue of $g^\ell$ at $x$}\Big\} \geq c_0>0.$$ Then, after passing to subsequences, $\left\{f^\ell\right\}$ converges weakly in $W^{2,p}$ to an immersion $\overline{f}: \M \emb \R^{n+k}$, whose induced metric $\overline{f}^\#\delta$, second fundamental form, and affine normal connection are limits of $\left\{g^\ell\right\}$,  $\left\{\two^\ell\right\}$, and  $\left\{\na^{\perp,\ell}\right\}$ in the weak $W^{1,p}$, $L^p$, and $L^p$-topologies, respectively. \end{corollary}

\begin{proof}
It follows immediately from the statement and the proof of  Theorem~\ref{thm: existence}.
\end{proof}

\subsection{A remark on smooth metrics without smooth isometric immersions}

In view of the smoothability Theorem~\ref{thm: main}, let us give a few tentative observations on the somewhat surprising examples of smooth metrics which only admit isometric immersions/embeddings of low regularity.

To put things into perspective, recall the following due to Burago--Shefel \cite{bs}; Iaia \cite{iaia}:
\begin{quote}
There exist real-analytic metrics on $\stwo$ with positive Gauss curvature except at one point, which admit only $C^{2,1}$ but no $C^3$-isometric embedding into $\R^3$.
\end{quote}
This leads to the question as follows, posed in \cite{gl} by Guan and Li:
\begin{quote}
Under what conditions on a smooth metric $g$ on $\stwo$ with nonnegative Gauss curvature, is there a $C^{2,\alpha}$ (even $C^{2,1}$) global isometric embedding into $\R^3$?
\end{quote}
Any obstruction to the existence of smooth isometric embeddings for nonnegatively curved metrics on surfaces must be global, for C.-S. Lin \cite{lin} established the local existence of smooth isometric embeddings in this case. The above examples of Pogorelov \cite{p} and Iaia \cite{iaia} further show that  such global obstruction cannot merely be $\pi_1(\M) \neq \{0\}$. 
On the other hand, example of a smooth, nonpositively curved metric on a surface that admits no \emph{local} isometric immersions into $\R^3$ was constructed by Khuri \cite{khuri}.

In this regard, we consider any given  $C^\infty$-Riemannian metric $g$ on $\M^n$ that admits a $W^{2,p}$- but no $C^\infty$-isometric immersion into $\R^{n+k}$; $p>n$.  We first make an almost tautological observation: for any such $W^{2,p}$-isometric immersion, there is no smoothing of the associated connection 1-forms which only modifies the extrinsic geometry while leaving invariant the intrinsic geometry. Equivalently:
\begin{quote}
Let $\M$ be a connected and simply-connected closed manifold of dimension $n$. Assume that $g$ is a $C^\infty$-Riemannian metric on $\M$ which admits a $W^{2,p}$- but no $C^\infty$-isometric immersion into $\R^{n+k}$; here $p>n$. Denote by $\Omega \in L^p(\M;T^*\M \otimes \mathfrak{so}(n+k,\R))$ the connection 1-form associated to the Cartan formalism of the $W^{2,p}$-isometric immersion  as in Equation~\eqref{Omega, new'}. Then, any smoothing sequence $\left\{\Omega^\e\right\}_{0<\e\leq 1} \subset C^\infty(\M;T^*\M \otimes \mathfrak{so}(n+k,\R))$ such that $\Omega^\e \to \Omega$ in $L^p$ on $\M$ (existence ensured by Theorem~\ref{thm: main}) cannot have the same upper-left block as $\Omega$. 
\end{quote}

Now let us re-examine the proof of Theorem~\ref{thm: main}. Recall Eq.~\eqref{Omega, new'} that 
\begin{equation*}
\Omega = \begin{bmatrix}
\na & \two\\
-\two & \na^\perp
\end{bmatrix} \in L^p_\loc\big(\M;T^*\M \otimes \mathfrak{so}(n+k,\R)\big),
\end{equation*}
where $\na$ and $\na^\perp$ are $n\times n$ and $k \times k$ blocks, respectively. (This is a schematic representation; see Eqs.~\eqref{Omega 1}--\eqref{Omega 3} for details. In particular, we refer to the $n\times n$ block $[\na]$ as the ``upper-left block''.) Also recall the Hodge decompositions in Step~1, Eq.~\eqref{new, decomposition of Omega} and Step~6, Eq.~\eqref{new: Omega epsilon def}:
\begin{align*}
\Omega = d\phi + d^*\psi + \harm,\qquad \Omega^\e = d\phi^\e + d^*\psi^\e + \harm,
\end{align*}
where $\phi^\e \to \phi$ in $W^{1,p}$ (Step~2) and $\psi^\e \to \psi$ in $W^{1,p}$ (Step~6); $\harm$ is a harmonic (hence $C^\infty$) 1-form. Here $\phi^\e$ are matrix-valued scalar fields, hence can be chosen to coincide with $\phi$ on the upper-left block. We assume that $\phi$ and $\psi$ are both $C^\infty$ in the upper-left block, for the same block of $\Omega$ (which contains components of the Levi-Civit\`{a} connection of $g$) is $C^\infty$.

However, the situation is different for the matrix-valued 2-forms $\psi^\e$: 
\begin{align}\label{word}
&\text{one cannot choose $\psi^\e$ such that their $n \times n$ upper-left blocks}\nonumber\\
&\qquad\qquad\text{coincide with that of $\psi$ for every $0<\e\leq 1$.}
\end{align}
This can be seen from Eq.~\eqref{system for psi epsilon}, the defining equation for $\psi^\e$ reproduced below\footnote{Recall here Remark~\ref{rem: dropping []}: the second term on the right-hand side is the intertwining of the wedge product on the differential form factor and Lie bracket on the matrix factor. }: 
\begin{equation*}
\Delta\psi^\e = \curv - \left(d\phi^\e + d^* \psi^\e + \harm\right) \wedge \left(d\phi^\e + d^* \psi^\e+\harm\right). 
\end{equation*}
By considering the $(i,j)$-entries,  $1 \leq i,j \leq n$,  on both sides of the above identity, we find that even the $(i,j)$-entry of $\Delta \psi^\e$ depends on the $(i, \alpha)$- and $(j,\beta)$-entries of $d^*\psi^\e$ for some $n+1 \leq \alpha, \beta \leq n+k$. These entries of $d^*\psi^\e$ are incorporated into the second fundamental forms of the smoothed immersions.

In summary, we may regard Eq.~\eqref{word} as a mechanism for the existence of $C^\infty$-Riemannian metric $g$ on $\M^n$ that admits a $W^{2,p}$- but no $C^\infty$-isometric immersion into $\R^{n+k}$ ($p>n$). It reflects the basic principle that the intrinsic and extrinsic geometries of isometric immersions are interrelated, which is also manifested in the structure of the Gauss--Codazzi--Ricci equations.

\subsection{Weak continuity of $W^{2,2+\e}$-isometric immersions revisited} 

In this subsection, we briefly recall the weak continuity theorem of the  Gauss--Codazzi--Ricci equations in \cite{csw}:
\begin{quote}
Weak $L^{2+\e}_\loc$-limits of (approximate) weak solutions to the Gauss--Codazzi--Ricci equations are also weak solutions, regardless of the dimension and codimension of the associated isometric immersions.
\end{quote}

This result has been revisited and extended in \cite{cl, l, giron} and other works. It is proved by way of utilising the div-curl lemma of Murat \cite{m} and Tartar \cite{t1, t2}, and it  serves as a key ingredient of the \emph{compensated compactness framework} for  the existence of isometric immersions with low regularity, pioneered by Chen--Slemrod--Wang \cite{csw'}: one constructs (\emph{e.g.}, by parabolic regularisation via adding artificial viscosity) approximate solutions to the Gauss--Codazzi--Ricci equations, obtains uniform $L^p_\loc$-estimates ($p>2$; \emph{e.g.}, via the method of invariant regions) for the approximate solutions, and then passes to the weak limit. In light of this weak continuity theorem \cite{csw}, any such weak limit constitutes a weak solution to the Gauss--Codazzi--Ricci equations.

The above result, nevertheless, is purely a PDE theorem: it only ascertains the weak continuity of solutions to the Gauss--Codazzi--Ricci equations, but cannot be translated into a geometric result, \emph{i.e.}, a weak stability statement for isometric immersions, unless $p>n$ as per Theorems~\ref{thm: existence, early} $\&$ \ref{thm: existence}. See Litzinger \cite{l} for the affirmative result in the endpoint case $p=n=2$.

In the spirit of Chen--Slemrod--Wang \cite{csw}, one may easily obtain a weak continuity theorem for the curvature Eq.~\eqref{curvature, def} for arbitrary curvature 2-form $\curv$. Again, the critical regularity index is $p=2$, regardless of the dimension of the manifold $\M$ or the rank of the bundle ${\rm Ad}\,\eta$.


\begin{theorem}\label{thm: weak continuity}
Let $\eta$ be a vector bundle with structure group $G$ over a Riemannian $n$-manifold $\M$, $n \geq 2$. 
Let $\{\Omega_j\}\subset L^p_\loc\left(\M;{\rm Ad}\,\eta\otimes \bigwedge^2 T^*\M\right)$ be a sequence of connection 1-forms satisfying
\begin{itemize}
\item
$\Omega_j \weak \overline{\Omega}$ weakly in the $L^p_\loc$-topology; and 
\item
$d\Omega_j + \Omega_j \wedge \Omega_j = \curv + \mathfrak{o}_j$ in the sense of distributions, where the curvature 2-form $\curv$ lies in $L^{p/2}_\loc$ and the error term $\mathfrak{o}_j \to 0$ in  $W^{-1,\sigma}$-topology for some $1 <\sigma \leq 2$ as $j \to \infty$.
\end{itemize}
Then we have $d\overline{\Omega}+\overline{\Omega}\wedge\overline{\Omega}=\curv$ in the sense of distributions.
\end{theorem}

\begin{proof}
It suffices to pass to the limit in the quadratic nonlinear term, \emph{i.e.}, to prove that
\begin{align*}
\Omega_j \wedge \Omega_j \, \longrightarrow \,
\overline{\Omega}\wedge\overline{\Omega}\quad\text{ in the sense of distributions}.
\end{align*}
In view of the wedge product theorem of compensated compactness (Robbin--Rogers--Temple \cite{rrt}, which is a generalisation of the classical div-curl lemma), it is enough to show that $\{d\Omega_j\}$ is precompact in the $W^{-1,2}_\loc$-topology.

To this end, we invoke the assumption $$d\Omega_j = -\Omega_j \wedge \Omega_j + \curv + \mathfrak{o}_j.$$ Fix any compact set $\K \Subset \M$. Then $\left\{\Omega_j \wedge \Omega_j + \curv\right\}$ is precompact in $L^{p/2}\left(\K; {\rm Ad}\,\eta\otimes \bigwedge^2 T^*\M\right)$ by the Cauchy--Schwarz inequality and the assumption on $\curv$. Thus, in view of the assumption on $\mathfrak{o}_j$ and the Sobolev embedding and/or the Rellich--Kondrachov lemma, one can find an index $\sigma' \in [1,2[$ such that $\{d\Omega_j\}$ is precompact in $W^{-1,\sigma'}\left(\K; {\rm Ad}\,\eta\otimes \bigwedge^2 T^*\M\right)$. On the other hand, since $\Omega_j \weak \overline{\Omega}$ weakly in $L^p\left(\K; {\rm Ad}\,\eta\otimes T^*\M\right)$ for  $p>2$, we infer that $\{d\Omega_j\}$ is bounded in $W^{-1,p}\left(\K; {\rm Ad}\,\eta\otimes \bigwedge^2 T^*\M\right)$. An interpolation theorem for negative-order Sobolev spaces (see, \emph{e.g.}, \cite[Theorem~3]{cdl} by Ding--Chen--Luo) shows that $\{d\Omega_j\}$ is precompact in $W^{-1,2}\left(\K; {\rm Ad}\,\eta\otimes \bigwedge^2 T^*\M\right)$.

Therefore, the assertion now follows from the wedge product theorem in \cite{rrt}.   \end{proof}

\section{Concluding remarks}\label{sec: remarks}


This note is largely motivated by S. Mardare \cite{m05, m07}, in which the fundamental theorem of surface theory (\emph{i.e.}, the equivalence between the existence of isometric immersions of connected, simply-connected closed surfaces into $\R^3$ and Gauss--Codazzi--Ricci equations) has been extended to isometric immersions with $W^{2,p}$-regularity; $p>2$. 
Recently, Litzinger \cite{l} further extended the fundamental theorem of surface theory  to $W^{2,2}$-isometric immersions, by way of exploiting Rivi\`{e}re's seminal work \cite{riviere} on the existence of nice gauges for conservation laws arising from geometric variational problems; see also Schikorra \cite{armin}.  Our note has taken a different approach: instead of using gauges, we perform a direct smoothing argument to Cartan's second structural equation~\eqref{curvature, def} --- $d\Omega + \Omega \wedge \Omega = \curv$ --- that \emph{fixes} $\curv$. This approach, unfortunately, only gives us existence results for $W^{2,p}$-isometric immersions $\M^n \emb \R^{n+k}$ in the \emph{subcritical} regime ($p>n$);  on the other hand, it nonetheless works for arbitrary dimension $n$ and codimension $k$.

It would be interesting to explore if the theories and methods in \cite{l, riviere, armin} may help tackle the \emph{critical} case $p=n$ for any $n$ and $k$. We hope to explore in future works whether gauge transforms \emph{\`{a} la} Uhlenbeck \cite{u', u} (see also \cite{w, w-new}) applied to Cartan's second structural equation~\eqref{curvature, def} may shed new lights on existence, rigidity, and other properties of isometric immersion with arbitrary dimension and codimension, even into general target manifolds. Investigations in this direction have recently been carried out in Giron's D.Phil. thesis \cite{giron} and Chen--Giron \cite{cg}.

It has also been brought to our attention that, in a series of recent  works \cite{rt1, rt2, rt3, rt4}, Reintjes--Temple systematically studied the question of smoothability of $L^\infty$-connections on tangent bundles of arbitrary Riemannian or pseudo-Riemannian manifolds with $L^\infty$-curvature tensor,  via a system of elliptic PDEs (``RT equations''). This furnishes a generalisation of Uhlenbeck compactness theorems \cite{u', u} and has novel applications to  general relativity. In particular, it is shown in \cite{rt2, rt3} that the Lorentzian metrics of shock wave solutions of the Einstein-Euler equations are non-singular. We wonder if the nonlinear smoothing techniques in \cite{m05, m07}, or the more geometrically-oriented results in our note, may find further applications in mathematical general relativity and/or other fields of physics.

\appendix
\section{Rudiments of differential geometry}\label{appendix}

For convenience of the reader, we briefly recall some notions and constructions in differential geometry used in this paper. More comprehensive treatment can be found in do Carmo~\cite{d} (on basics of Riemannian geometry) and Wehrheim~\cite{w} (on principal $G$-bundle and Coulomb--Uhlenbeck gauge), among many other excellent texts.

An $n$-dimensional Riemannian manifold $(\mathcal{M},g)$ of regularity $\mathcal{X}$ is a second-countable, Hausdorff topological space with an atlas of local charts $\mathcal{A} = \{(U_{\alpha},\phi_{\alpha}):\alpha \in \mathcal{I}\}$ such that each $U_{\alpha}\subset \mathcal{M}$ is an open subset, $\phi_{\alpha}:U_{\alpha}\to \phi_{\alpha}(U_{\alpha})\subset \mathbb{R}^{n}$ is a homeomorphism, and the transition functions $\phi_{\alpha}\circ \phi_{\beta}^{-1}$ between the overlapping charts $U_{\alpha}$ and $U_{\beta}$ lie in the regularity class $\mathcal{X}$. For instance, we may consider $\mathcal{X}=C^\infty, C^2, C^1, C^0, W^{2,p},\ldots$. If each $\phi_{\alpha}\circ \phi_{\beta}^{-1}$ can be chosen to have positive definite Jacobian determinant almost everywhere, then $\mathcal{M}$ is said to be orientable. The Riemannian metric $g$ on $\mathcal{M}$ is given by a field of inner products --- at each point $P\in \mathcal{M}$, $g(P):\Gamma (T\mathcal{M})\times \Gamma (T\mathcal{M})\to \mathbb{R}$ is an inner product denoted by
\[
g(P)(X,Y) = \langle X,Y\rangle \qquad \mathrm{for~any~}X,Y\in \Gamma (T\mathcal{M}),
\] 

Throughout, $\Gamma (T\mathcal{M})$ denotes the space of vector fields on $\mathcal{M}$. More generally, we denote by $\Gamma$ the space of sections of given vector or fibre bundles with required (Sobolev) regularity; here $T\M$ is the tangent bundle of the manifold $\M$.

Given $(\mathcal{M},g)$, an affine connection on $T\mathcal{M}$ is a bilinear map $\nabla :\Gamma (T\mathcal{M})\times \Gamma (T\mathcal{M})\to \Gamma (T\mathcal{M})$ such that, for any $f:\mathcal{M}\to \mathbb{R}$ and any $X,Y,Z\in \Gamma (T\mathcal{M})$,
\[
\nabla_{fX}Y = f\nabla_{X}Y,\qquad \nabla_{X}(fY) = f\nabla_{X}Y + X(f)Y,
\]
where $X(f)$ is the directional derivative of $f$ in the direction of $X$. We say that $\nabla$ is compatible with metric $g$ if
\[
Xg(Y,Z) = g(\nabla_{X}Y,Z) + g(Y,\nabla_{X}Z),
\]
and $\nabla$ is torsion-free if
\[
\nabla_{X}Y - \nabla_{Y}X = [X,Y],
\]
where the Lie bracket is defined by $[X,Y] = XY - YX$. There exists a unique compatible, torsion-free affine connection $\nabla$ on $\mathcal{M}$ known as the Levi-Civita connection, where the bilinear map $\nabla$ is also called the covariant derivative. As a basic example, consider $(\bar{M},\bar{g}) = (\mathbb{R}^{m},g_{0})$ with the Euclidean metric $g_{0} = \delta_{ij}$, i.e., the dot product, whose Levi-Civita connection $\bar{\nabla}$ is given by $\bar{\nabla}_{X}Y:= XY$.

Throughout this work, as is customary, we always identify a vector field $X \in \G(T\M)$ as a field of first-order differential operator. That is, for any smooth function (scalar field) $f: \M \to \R$, we have $Xf \equiv X^i\p_if$ for $X=X^i\p_i$ in any local coordinate system $\{\p_i\}_1^n \subset \G(T\M)$. 
 
Given a manifold $\mathcal{M}$, we say that $(E,\mathcal{M},F)$ is a vector bundle of degree $k\in \mathbb{N}$ over $\mathcal{M}$ if there is a surjection $\pi :E\to \mathcal{M}$ such that, for any $P\in \mathcal{M}$, there exists a local neighborhood $U\subset \mathcal{M}$ containing $P$ so that there is a diffeomorphism $\psi_{U}: \pi^{-1}(U)\to U\times F$ with $\mathrm{pr}_{1}\circ \psi_{U} = \pi$ on $\pi^{-1}(U)$, where map $\mathrm{pr}_{1}$ is the projection map onto the first coordinate, $E$ and $F$ are differentiable manifolds, and $F \simeq \mathbb{R}^{k}$ as a vector space isomorphism. In this bundle, $E$ is called the total space, $F$ is the fibre, $\mathcal{M}$ the base manifold, $\pi$ the projection of the bundle, and $\psi_{U}$ a local trivialisation. For simplicity, we also say that $E$ is a vector bundle over $\mathcal{M}$.

If $E_{1}$ and $E_{2}$ are both vector bundles over $\mathcal{M}$, we can take the direct sum and the quotient of the bundles, by taking the vector space direct sum and quotient of the fibres. Also, for a vector bundle $(E,\mathcal{M},F = \mathbb{R}^{k})$ with projection $\pi$, the space of smooth sections is defined by $\Gamma (E) := \{s \in C^{\infty}(\mathcal{M}; E) : \pi \circ s = \mathrm{id}_{\mathcal{M}}\}$. We can define the affine connection $\nabla^{E} : \Gamma (T\mathcal{M}) \times \Gamma (E) \mapsto \Gamma (E)$, by linearity and the Leibniz rule. 

As a primary example, consider $E = T\mathcal{M}$, the tangent bundle over $\mathcal{M}$. Then $\pi$ is the projection onto the base point in $\mathcal{M}$, and $\Gamma (T\mathcal{M})$ is the space of smooth vector fields. Moreover, $\nabla^{T\mathcal{M}}$ is precisely the Levi-Civita connection. Another example is the cotangent bundle $T^{*}\mathcal{M}$ over $\mathcal{M}$, whose fibres are the dual vector spaces of the fibres of $T\mathcal{M}$.

As another important example, consider $E_{1} = T^{*}\mathcal{M} \otimes T^{*}\mathcal{M} \otimes \dots \otimes T^{*}\mathcal{M}$, the tensor product of $q$-copies of $T^{*}\mathcal{M}$, for $q = 0, 1, 2, \ldots$. This is the (covariant) $q$-tensor algebra over $\mathcal{M}$, which can be viewed as $q$-linear maps on $T\mathcal{M}$. Now let $E_{2} \subset E_{1}$ be the subspace of all the alternating $q$-tensors on $T\mathcal{M}$, i.e., the $q$-linear maps that change sign when switching any pair of its indices $\{i, j\} \subset \{1, \ldots , q\}$. By convention, we write $E_{2} := \bigwedge^{q} T^{*}\mathcal{M}$, known as the alternating $q$-algebra. Moreover, the sections are $\Omega^{q}(\mathcal{M}) := \Gamma (\bigwedge^{q} T^{*}\mathcal{M})$, known as the differential $q$-forms. A generic element $\alpha \in \Omega^{q}(\mathcal{M})$ is a linear combination of alternating forms $\xi_{1} \wedge \dots \wedge \xi_{q}$, where $\{\xi_{j}\}_{j = 1}^{q}$ is a $q$-tuple of linearly independent differential 1-forms on $\mathcal{M}$. 

On the space of differential forms, we recall four important operations. The first is the exterior derivative $d: \Omega^{q}(\mathcal{M}) \to \Omega^{q + 1}(\mathcal{M})$, which is linear and satisfies $d^{2} = 0$. The second is the Hodge star $*: \bigwedge^{q} T^{*}\mathcal{M} \to \bigwedge^{n - q} T^{*}\mathcal{M}$, which can also be regarded as $*: \Omega^{q}(\mathcal{M}) \to \Omega^{n - q}(\mathcal{M})$. It is an isomorphism of vector bundles, which satisfies $** = (- 1)^{q(n - q)}$ whenever $\mathcal{M}$ is orientable. The third is a natural product on differential forms: For $\alpha \in \Omega^{q}(\mathcal{M})$ and $\beta \in \Omega^{r}(\mathcal{M})$, we can define the wedge product $\alpha \wedge \beta \in \Omega^{q + r}(\mathcal{M})$. The fourth is the covariant derivative $\nabla : \Omega^{q}(\mathcal{M}) \to \Omega^{q + 1}(\mathcal{M})$: For $\alpha \in \Omega^{q}(\mathcal{M})$, define $\nabla \alpha \in \Omega^{q + 1}(\mathcal{M})$ via the Leibniz rule:
\begin{align*}
(\nabla \alpha)(X,Y_{1},\ldots ,Y_{q}) &\equiv \nabla_{X}\alpha (Y_{1},\ldots ,Y_{q})\\
&:= X\big(\alpha (Y_{1},\ldots ,Y_{q})\big) - \alpha (\nabla_{X}Y_{1},\ldots ,Y_{q}) - \dots - \alpha (Y_{1},\ldots ,\nabla_{X}Y_{q})
\end{align*}
for any $X, Y_{1}, \ldots , Y_{q} \in \Gamma (T\mathcal{M})$.

There is a natural isomorphism between $T\mathcal{M}$ and $T^{*}\mathcal{M}$, by identifying canonically each fibre of $T\mathcal{M}$ with its dual. It induces a canonical isomorphism $\sharp \equiv \flat : \Omega^{1}(\mathcal{M}) \to \Gamma (T\mathcal{M})$. When a vector field and its corresponding 1-form in local coordinates are written by the Einstein summation convention, $\sharp$ (or $\flat$) amounts to raising (or lowering) the indices of the coefficients. Clearly, they extend to the isomorphisms between $T^{*}\mathcal{M}\otimes \dots \otimes T^{*}\mathcal{M}$ (covariant tensors) and $T\mathcal{M}\otimes \dots \otimes T\mathcal{M}$ (contravariant tensors).

For instance, consider the covariant derivative $\nabla :\Omega^{q}(\mathcal{M})\to \Omega^{q + 1}(\mathcal{M})$ defined above. For $q = 1$ and $\alpha \in \Omega^{1}(\mathcal{M})$, we set $X:= \alpha^{\sharp}\in \Gamma (T\mathcal{M})$. Since $\nabla_{Y}X\in \Gamma (T\mathcal{M})$ for any $Y\in \Gamma (T\mathcal{M})$, we can view $\nabla X = \nabla \alpha^{\sharp}\in \Gamma (T\mathcal{M}\otimes T\mathcal{M})$, i.e., $(\nabla \alpha^{\sharp})^{\flat}\in \Omega^{2}(\mathcal{M})$. This example shows that, via the identifications $\sharp$ and $\flat$, the covariant derivative $\nabla$ on $\Omega^{q}(\mathcal{M})$ generalizes the definition of the Levi-Civita connection.

On an orientable $n$-dimensional Riemannian manifold $(\mathcal{M},g)$, there is a natural $n$-form $\mathrm{d}V_{g}\in \Omega^{n}(\mathcal{M})$, called the Riemannian volume form, which satisfies $*1 = \mathrm{d}V_{g}$ and $*\mathrm{d}V_{g} = 1$. Let $\mathcal{A} = \{(U_{\alpha},\phi_{\alpha}):\alpha \in \mathcal{I}\}$ be an oriented atlas for $\mathcal{M}$. There exists a locally finite $C^{\infty}$-partition of unity $\{\rho_{\alpha}:\alpha \in \mathcal{I}\}$  subordinate to $\mathcal{A}$; namely that $\sum_{\alpha \in \mathcal{I}}\rho_{\alpha} = 1$, $0\leq \rho_{\alpha}\leq 1$, and $\mathrm{supp}(\rho_{\alpha})\Subset U_{\alpha}$. We define the integration of $\omega \in \Omega^{n}(\mathcal{M})$ over $\mathcal{M}$ by
\[
\int_{\mathcal{M}}\omega := \sum_{\alpha \in \mathcal{I}}\int_{\mathbb{R}^{n}}\rho_{\alpha}((\phi_{\alpha}^{-1})^{\#}\omega)\,\chi_{\phi_{\alpha}(U_{\alpha})}\sqrt{|g|}\,\mathrm{d}x_{1}\cdots \mathrm{d}x_{n},
\]
where $\phi^{\#}$ denotes the pullback of a tensor under $\phi$ on $\mathcal{M}$ and $|g| := \det(g_{ij})$. The integration of function $\phi$ on $\mathcal{M}$ is defined as the integration of its Hodge dual, i.e., $\int_{\mathcal{M}}\phi := \int_{\mathcal{M}}\phi \,\mathrm{d}V_{g}$.

We now define the Sobolev spaces $W^{k,p}(\mathcal{M}; \bigwedge^{q}T^{*}\mathcal{M})$ on an $n$-dimensional Riemannian manifold $(\mathcal{M},g)$, which generalise $W^{k,p}(\mathbb{R}^{n})$ for $k \in \mathbb{Z}$ and $1 \leq p \leq \infty$. First, for differential $q$-forms $\alpha , \beta \in \Omega^{q}(\mathcal{M})$, $g$ on $\mathcal{M}$ defines an inner product $g(\alpha , \beta) = \langle \alpha , \beta \rangle$ by
\[
\langle \alpha ,\beta \rangle \,\mathrm{d}V_{g}:= \alpha \wedge *\beta.
\]
Then, for alternating contravariant $q$-tensor fields $T$ and $S$, we set $\langle T,S\rangle := \langle T^{\flat},S^{\flat}\rangle$. Next, the $L^{p}$-norm of $\alpha \in \Omega^{q}(\mathcal{M})$ is defined as
\[
\| \alpha \|_{L^{p}} := \left(\int_{\mathcal{M}}\left[*\left(\alpha \wedge *\alpha\right)\right]^{\frac{p}{2}}\mathrm{d}V_{g}\right)^{\frac{1}{p}} = \left(\int_{\mathcal{M}}\langle \alpha ,\alpha \rangle^{\frac{p}{2}}\mathrm{d}V_{g}\right)^{\frac{1}{p}}.
\]
Moreover, for $\alpha \in \Omega^{q}(\mathcal{M})$, set
\[
\| \alpha \|_{W^{k,p}} := \left(\sum_{j = 0}^{k}\| \nabla^{j}\alpha \|_{L^{p}}^{p}\right)^{\frac{1}{p}} = \left(\sum_{j = 0}^{k}\| \underbrace{\nabla \circ \dots \circ \nabla}_{j\text{ times}}\alpha \|_{L^{p}}^{p}\right)^{\frac{1}{p}}.
\]
We denote by $W^{k,p}(\mathcal{M}; \bigwedge^{q}T^{*}\mathcal{M})$ the completion of the space of compactly supported $q$-forms with respect to $\| \cdot \|_{W^{k,p}}$. Notice that
\[
\| X\|_{W^{k,p}} := \| X^{\flat}\|_{W^{k,p}}
\]
for any contravariant tensor field $X$. In fact, $W^{k,p}(\mathcal{M}; E)$ can be defined for an arbitrary bundle $E$. Furthermore, we can define (indeed, intrinsically) the $W^{k,p}$-connections $\nabla^{E}$ on any vector or fibre bundle $E$. This can be done since the moduli space of connections is an affine space modelled over the tensor algebra $\bigotimes^{\bullet} T\mathcal{M} \otimes \bigotimes^{\bullet} T^{*}\mathcal{M} \otimes E$. 

The construction of vector bundles can be further generalised in the following manner. Given a vector bundle $(E,\pi, \M; \R^k)$, one may replace the fibre (\emph{i.e.}, the vector space $\R^k$, which is an abelian Lie group) with a general Lie group $G$ that has a smooth right-action on $E$.  We again require $\pi:E \to \M$ to be a smooth surjection, and that for each $x \in \M$ there exist a local neighbourhood $U \subset \M$ containing $x$ and a diffeomorphism $\phi_U: \pi^{-1}(U) \to U \times G$. Furthermore, assume that for each $x \in \M$ the group $G$ acts freely and transitively on the fibre $\pi^{-1}\{x\}$, and the local trivialisation $\phi_U$ commutes with the $G$-action. In this way, $(E,\pi,M;G)$ is said to be a \emph{principal $G$-bundle} over the base manifold $\M$ with projection $\pi$ and \emph{structure group} $G$.  Following the  seminal works \cite{u, u'} by Uhlenbeck, we write $\eta$ for $E$ in this setting.

For a Lie group $G$, its Lie algebra $\mathfrak{g}$ is the tangent space at the identity element, encoding all ``infinitesimal directions'' of the group. The group $G$ acts on its own Lie algebra by \emph{conjugation}:
\[
\mathrm{Ad}: G \longrightarrow \mathrm{Aut}(\mathfrak{g}), \qquad \mathrm{Ad}_g(X) := gXg^{-1}.
\]
This action is called the \emph{adjoint representation} of $G$. Given a principal $G$-bundle $\pi: \eta \to M$, and given any space $F$ on which $G$ acts, one can construct an associated bundle with fibre $F$. In particular, if we take $F = \mathfrak{g}$ and let $G$ act on its Lie algebra via the adjoint representation, the corresponding associated bundle is called the \emph{adjoint bundle}:
\[
\mathrm{Ad} \,\eta := \eta \times_{\mathrm{Ad}} \mathfrak{g}.
\]
Its elements are equivalence classes $[p, X]$, where $p \in \eta$ and $X \in \mathfrak{g}$. Since the fibres inherit the Lie algebra structure, $\mathrm{Ad}\, \eta$ is a vector bundle whose fibres are Lie algebras.

\bigskip
\noindent
{\bf Acknowledgement}. SL thanks Armin Schikorra for many kind discussions on isometric immersions and gauge equations over the years, and Qing Han for insightful remarks on the existence theory of isometric immersions.

The research of SL is supported by NSFC Projects 12331008 $\&$ 12411530065, the Young Elite Scientists Sponsorship Program by CAST 2023QNRC001, National Key Research $\&$ Development Programs 2023YFA1010900 and 2024YFA1014900, Shanghai Rising-Star Program 24QA2703600, Shanghai Qi-Guang Scholarship, and Shanghai Frontiers Science Center of Modern Analysis.

\medskip
\noindent
{\bf Statement of competing interests}. We declare that there are no conflicts of interest involved.

\medskip
\noindent
{\bf Data Availability Statement}. We declare that no data are associated with this work.

\medskip
\noindent
{\bf AI Statement}. The author confirms that no generative AI or AI-assisted technologies were used in the writing or preparation of this manuscript.

\end{document}